\documentclass[twocolumn]{svjour3}          
\smartqed  
\usepackage{graphicx}

\usepackage{latexsym,amsmath,amssymb,amsfonts,graphics,epsfig}
\usepackage{float}
\usepackage{pdflscape}
\usepackage{afterpage}
\usepackage{epstopdf}

\begin{document}

\title{Analysis of Neural Clusters due to \\ Deep Brain Stimulation Pulses\\}

\author{Daniel Kuelbs \and Jacob Dunefsky \and Bharat Monga \and Jeff Moehlis}

\institute{
Daniel Kuelbs \at Stanford University, Palo Alto, CA, 94305 \\ \email{dkuelbs@stanford.edu}, Co-first Author
\and
Jacob Dunefsky \at Yale University, New Haven, CT, 06520 \\ \email{jacob.dunefsky@yale.edu}, Co-first Author
\and
Bharat Monga \at Department of Mechanical Engineering, University of California, Santa Barbara, CA 93106 \\ \email{monga@ucsb.edu}
\and
Jeff Moehlis \at Department of Mechanical Engineering, Program in Dynamical Neuroscience, University of California, Santa Barbara, CA 93106 \\ \email{moehlis@ucsb.edu}, Corresponding Author
}

\maketitle


\begin{abstract}

Deep brain stimulation (DBS) is an established method for treating
pathological conditions such as Parkinson's disease, dystonia, Tourette 
syndrome, and essential tremor.  While the precise mechanisms which underly
the effectiveness of DBS are not fully understood, theoretical studies of
populations of neural oscillators stimulated by periodic pulses suggest that 
this may be related to clustering, in which subpopulations of the neurons 
are synchronized, but the subpopulations are desynchronized with respect to 
each other.  The details of the clustering behavior depend on the frequency 
and amplitude of the stimulation in a complicated way.  In the present study,
we investigate how the number of clusters, their stability properties, and 
their basins of attraction can be understood in terms of one-dimensional maps 
defined on the circle.  Moreover, we generalize this analysis to stimuli that
consist of pulses with alternating properties, which provide additional 
degrees of freedom in the design of DBS stimuli.  Our results illustrate how
the complicated properties of clustering behavior for periodically forced 
neural oscillator populations can be understood in terms of a much simpler 
dynamical system.



\end{abstract}


\section{Introduction}

A primary motivation for this study is Parkinson's disease, which 
can cause an involuntary shaking that typically affects the distal 
portion of the upper limbs, and difficulty initiating motion.
For patients with advanced Parkinson's disease who do not respond to
drug therapy, electrical deep brain stimulation (DBS), an FDA-approved
therapeutic procedure, may offer relief \cite{bena91}.  Here, a neurosurgeon 
guides a small electrode into the sub-thalamic nucleus or globus pallidus 
interna (GPi); the electrode is connected to a pacemaker implanted 
in the chest which sends periodic electrical pulses directly into the brain 
tissue.  The efficacy of DBS for the treatment of Parkinson's disease has 
been found to depend on the frequency of stimulation, with high-frequency 
stimulation (70 to 1000 Hz and beyond) being therapeutically 
effective \cite{bena91,rizz01,moro02}. The generally accepted
therapeutic range is 130-180 Hz \cite{volk02,kunc04}.  
Experimental evidence has suggested that motor symptoms of Parkinson's 
disease are associated with pathological synchronization of neurons in 
the basal ganglia, and that DBS desynchronizes the neural 
activity~\cite{uhlh06,chen07,hamm07,levy00,schn05}.   DBS has also 
shown promising results in treating other neurological conditions,
for which the stimulation electrode is implanted in the GPi (for dystonia)
or the thalamus (for Tourette syndrome and essential 
tremor)~\cite{savi12,bena02}.


While the precise mechanisms which underly the effectiveness of DBS are not 
fully understood, theoretical studies have shown that DBS-like stimulation 
consisting of a periodic pulses applied to neural
oscillator populations can lead to chaotic desynchronization~\cite{wils11} 
or clustering behavior~\cite{wils15cluster}, in which subpopulations of 
the neurons are synchronized, but the subpopulations are desynchronized with 
respect to each other.  Clustering has also been found in theoretical studies 
of coordinated reset, in which multiple electrodes deliver inputs which are 
separated by a time delay~\cite{luck13,lysy11,lysy13,tass03a}.  These
studies, along with clinical successes with coordinated reset~\cite{adam14}, 
point to clustering as an attractive objective for designing stimulation 
properties; this has motivated the design of single control inputs which 
promote clustering~\cite{matc18,mong19_physicad,wils20}, in contrast to methods which 
seek to fully desynchronize the neural 
activity~\cite{tass03,nabi13,wils14a,mong20}.
Notably, clustering has at least two important differences from chaotic 
desynchronization: clustered states often exist over a much larger parameter 
range than chaotic desynchronization, a possible explanation why effective 
DBS parameters are easier to find than chaotic desynchronization would 
suggest; and clustered states may induce plasticity changes more effectively 
than chaotic desynchronization, which may explain why benefits are more 
persistent for some kinds of stimulation mechanisms than others 
(cf.~\cite{adam14,mong19_physicad}).  In this paper, we will focus on
clustering which arises from a single stimulation electrode, unlike
coordinated reset which uses multiple electrodes.

Despite substantial data backing the general efficacy of DBS, it can have
side effects including disorientation, memory deficits, spatial delayed 
recall, response inhibition, episodes of mania, hallucinations, or mood 
swings, as well as impairment of social functions such as the 
ability to recognize the emotional tone of a face~\cite{cyro16,buhm17}.  
Our study develops tools which can help to identify different stimuli that 
result in the same clustering behavior; our hope is that
the identification of these alternatives will allow neurologists to consider
different stimuli in order to find those which are effective at treating
neurological disorders while minimizing 
the severity of side effects.

In this paper, we investigate how the details of clustering due to periodic 
pulses of the type used in DBS can be understood in terms of 
one-dimensional maps defined
on the circle.  As a first step, Section~\ref{section:phase_reduction} 
describes phase reduction, a powerful classical technique for the analysis 
of oscillators in which a single variable describes the phase of the 
oscillation with respect to some reference state.  
Section~\ref{section:simulations} shows results from simulations of 
populations of neural ocsillators stimulated by periodic pulses
of the type used for DBS; this illustrates the different types of clustering 
which can occur, and motivates the theoretical analysis. 
Section~\ref{section:identical_stimuli} derives and investigates the 
one-dimensional maps which can be used to understand the types of clusters 
which occur, their stability properties, and their basins of attraction.  
Section~\ref{section:alternating_stimuli} then demonstrates how this 
analysis in terms of maps can be generalized to consider stimuli that
consist of pulses with alternating properties, which provide additional 
degrees of freedom for DBS stimulus design.
Section~\ref{section:conclusion} summarizes the results.  The models for 
the neurons considered in this paper are given in the Appendix.

\section{Phase Reduction} \label{section:phase_reduction}
A common way to describe the dynamics of neurons is to use conductance-based models such as the Hodgkin-Huxley equations~\cite{hodg52d}.  Such models are typically high-dimensional and contain a large number of parameters, which can make them unwieldy for simulations of large neural populations.  A powerful technique for the analysis of oscillatory neurons, whose dynamics are described by a stable periodic orbit, is the rigorous reduction of conductance-based models to phase models, with a single variable $\theta$ describing the phase of the oscillation with respect to some reference state~\cite{winf01,kura84,mong19}.

Suppose that our conductance-based model is described by the $n$-dimensional dynamical system 
\begin{equation}
\label{eq:Phase_reduction_1}
\frac{d {\bf x}}{dt}= {\bf F}({\bf x}) , \qquad {\bf x} \in \mathbb{R}^n \quad (n \geq 2),
\end{equation}
with a stable periodic orbit $\gamma(t)$ with period $T$.  For each point ${\bf x}^\ast$ in the basin of attraction of $\gamma(t)$ there exists a corresponding phase $\theta$(${\bf x}^\ast$) such that~\cite{guck75,winf01}
\begin{equation}
\label{eq:Phase_reduction_2}
\lim\limits_{t \to \infty} \left|{\bf x}(t) - \gamma \left( t+\frac{T}{2\pi}\theta({\bf x}^\ast) \right) \right| = 0,
\end{equation}
where, under the given vector field, ${\bf x}(t)$ is the trajectory of the initial point ${\bf x}^\ast$. The asymptotic phase of ${\bf x}$, $\theta({\bf x})$, ranges in value from $[0,2\pi)$. In this paper, $\theta = 0$ will represent the phase at which the neuron fires an action potential. Isochrons are level sets of $\theta({\bf x})$, and we define isochrons such that the phase of a trajectory evolves linearly in time both on and off of the periodic orbit~\cite{winf67,winf01}. As a result, for the entire basin of attraction of the periodic orbit,
\begin{equation}
\label{eq:Phase_reduction_3}
\frac{d\theta}{dt}= \frac{2\pi}{T} \equiv \omega.
\end{equation}

If we now consider the dynamical system
\begin{equation}
\label{eq:Phase_reduction_4}
\frac{d{\bf x}}{dt}= {\bf F}({\bf x}) + {\bf U}(t), \qquad {\bf x} \in \mathbb{R}^n,
\end{equation}
where ${\bf U}(t) \in \mathbb{R}^n$  is an infinitesimal control input, phase reduction gives the one-dimensional system~\cite{kura84,brow04,mong19}
\begin{equation}
\label{eq:Phase_reduction_5}
\frac{d\theta}{dt}= \omega + {\bf U}(t)^T {\bf Z}(\theta).
\end{equation}
In this equation, ${\bf Z}(\theta)$ is the gradient of $\theta$ evaluated on the periodic orbit, and is known as the phase response curve (PRC)~\cite{winf01,erme10,neto12}; it represents the change in phase that the control input will cause when applied at a given phase.  In this paper, we consider electrical current inputs which only act in the voltage direction defined by the unit vector $\hat{V}$, i.e., ${\bf U}(t) = u(t) \hat{V}$, with the corresponding phase reduction 
\begin{equation}
\label{eq:Phase_reduction_6}
\frac{d\theta_{i}}{dt}= \omega + Z(\theta_{i}) u(t).
\end{equation}
Here, $\theta_i$ represents the phase of the $i^{\rm th}$ neuron, $\omega$ is the natural frequency of the neuron in radians per second, $Z(\theta) = \frac{\partial \theta}{\partial V}$ is the component of the PRC in the voltage direction, and $u(t)$ is the input.  For the populations of neurons considered in this paper, we assume that the neurons are identical and they all receive the same input, and we will consider uncoupled neurons without noise; these assumptions allow a more detailed analysis to be performed.

In the next section, we show simulation results for populations of neurons described by such phase models with periodic pulses of the type used for DBS.

\section{Simulation Results for Identical Periodic DBS Pulses} \label{section:simulations}

\begin{figure}[tb]
\begin{center}
\leavevmode
\epsfxsize=3.2in
\epsfbox{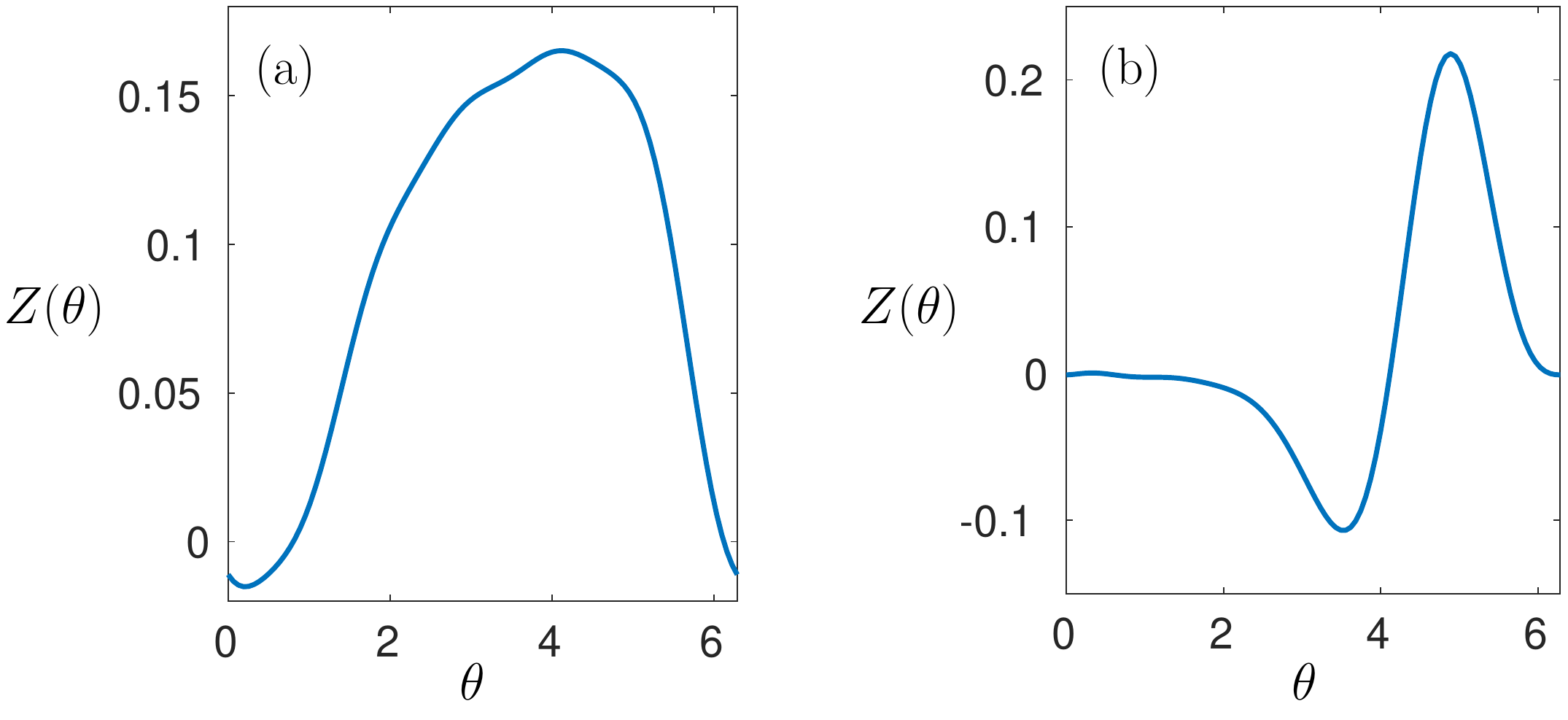}
\end{center}
\caption{Panels (a) and (b) show the phase response curves $Z(\theta)$ of the thalamic (Type I) and Hodgkin-Huxley (Type II) neurons considered in this paper, respectively. \label{Z}}
\end{figure}

In this section, we show simulation results for populations of neurons stimulated by periodic pulses of the type used for DBS; these results will inspire the analysis in Section~\ref{section:identical_stimuli}.  To illustrate a range of clustering behaviors, we show simulations for prototypical systems which represent two common types of neurons~\cite{rinz98}: as a Type I neuron model, we consider the model for thalamic neurons from~\cite{rubi04}, and as a Type II model we consider the Hodgkin-Huxley equations~\cite{hodg52d}.  These models are not meant to correspond to the neurons directly relevant to Parkinson's disease in human patients; rather, they are used to illustrate typical clustering behaviors for populations of neural oscillators under DBS-like stimuli.  The full equations and parameters for these models are given in the Appendix.  For our simulations, we use the corresponding phase models.  For reference, for these parameters the thalamic neurons have $\omega = 0.748$ rad/s, and the Hodgkin-Huxley neurons have $\omega = 0.429$ rad/s.  The PRC functions for these neurons are shown in Figure~\ref{Z}(a) and (b), respectively.  Each PRC was calculated numerically using XPP~\cite{erme02}, and is approximated by a Fourier series.

\begin{figure}[tb]
\begin{center}
\leavevmode
\epsfxsize=3.2in
\epsfbox{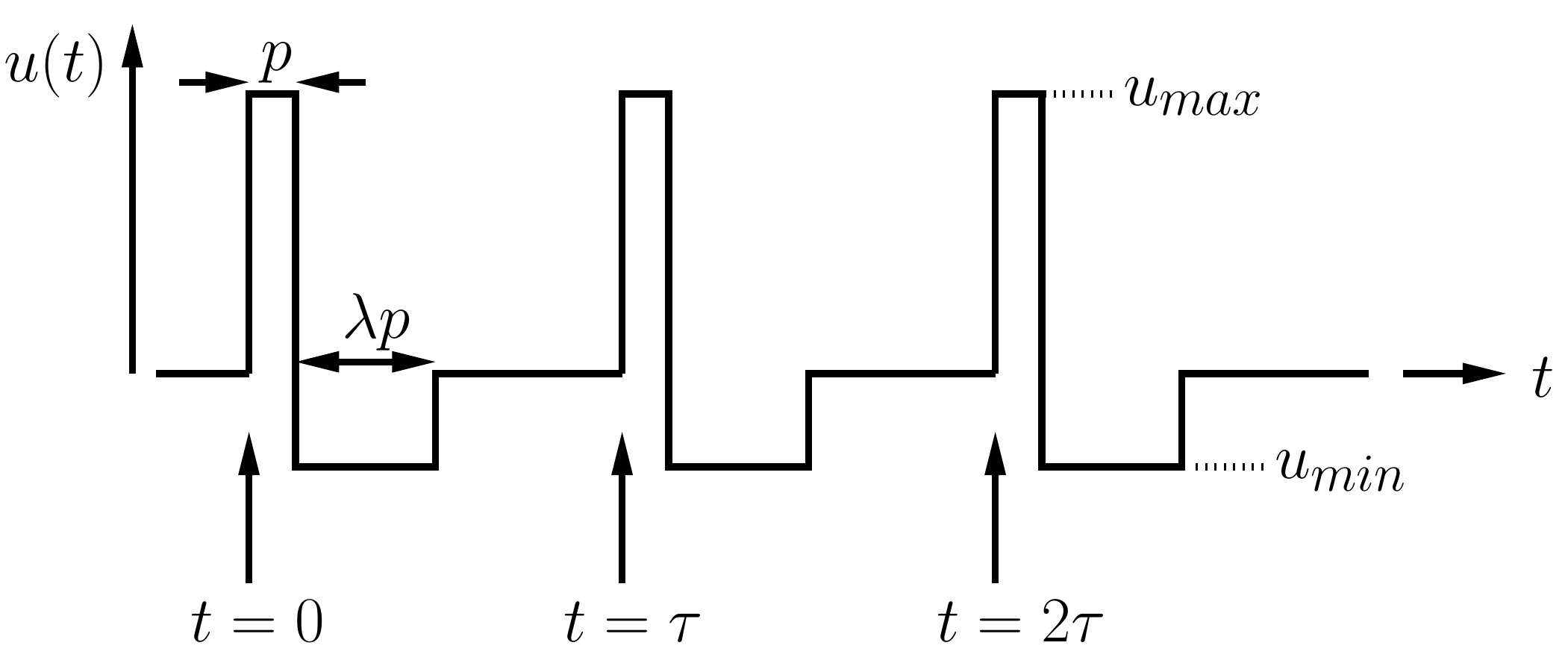}
\end{center}
\caption{Periodic sequence of identical pulses.
\label{one_kick}}
\end{figure}

The input $u(t)$ that we consider, shown in Figure~\ref{one_kick} and inspired by DBS stimuli~\cite{mont10}, is a periodic sequence of identical charge-balanced pulses parameterized by amplitude $u_{max}$, period $\tau$ (with corresponding frequency $1/\tau$), pulse width $p$, and multiplier $\lambda$ (the ratio of time that the pulse is negative to the time that the pulse is positive).  Mathematically, $u(t)$ is given by:
\begin{equation}
\label{eq:Pulsatile_Stimulus}
u(t) = \left\{
\begin{array}{ll}
u_{max} &  \bmod(t,\tau) \leq  p \\
u_{min} \equiv -\frac{u_{max}} {\lambda} & p < \bmod(t,\tau) \leq (\lambda+1)p\\
0 & $otherwise$.
\end{array}
\right.
\end{equation}
Unless otherwise stated, we will use $u_{max}$ corresponding to a current density of $20 \mu A/cm^2$, $p = 0.5$ ms, and $\lambda = 3$ in our simulations.  We consider different frequencies of stimulation between 70-300 Hz, which includes the typical therapeutic range of 130-180 Hz for DBS treatment of Parkinson's disease.

\begin{figure}[!t]
\begin{center}
\includegraphics[width=0.45\textwidth]{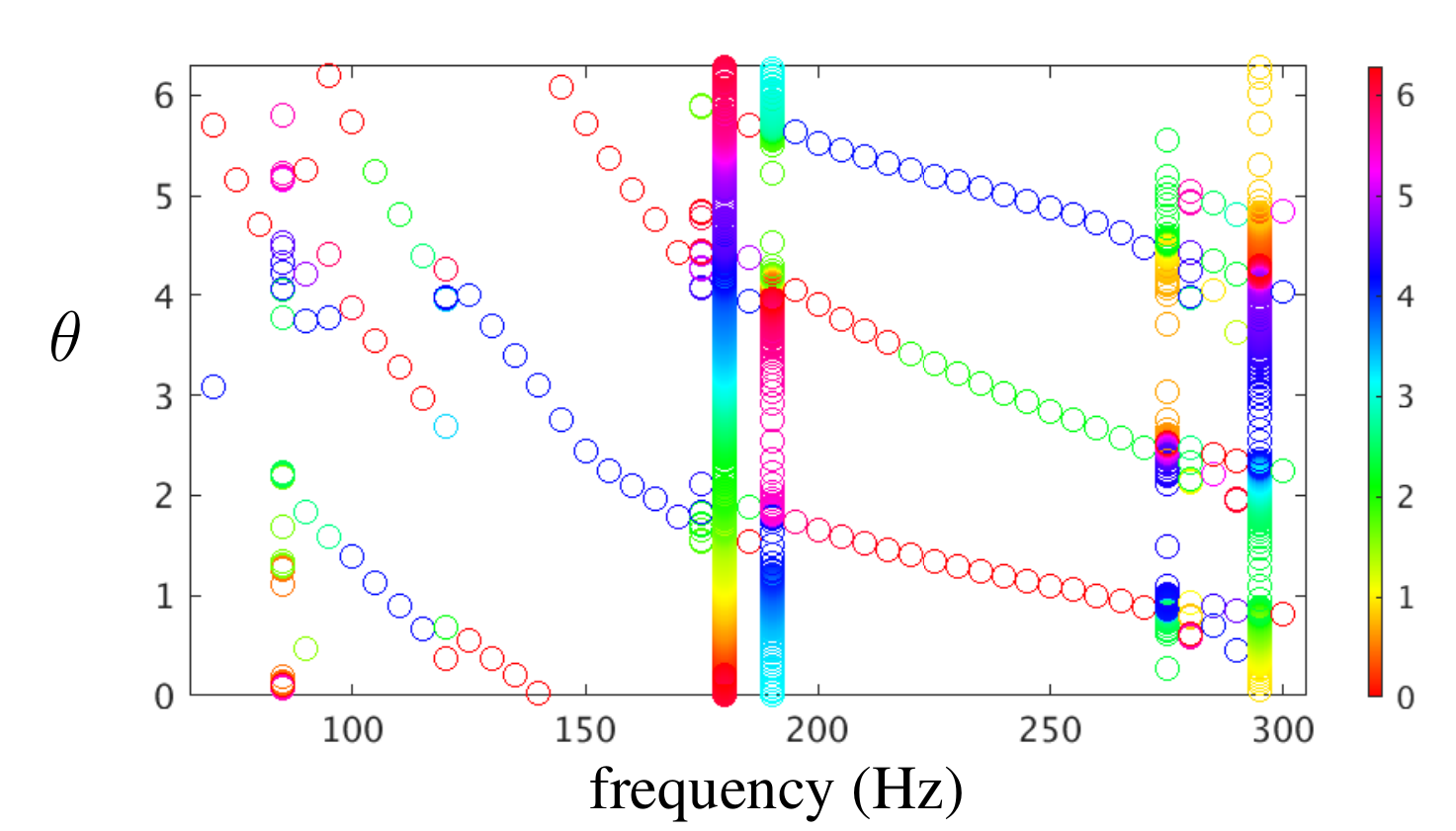}
\end{center}
\caption{The final phases $\theta$ of Hodgkin-Huxley neurons drawn from an initial uniform distribution as a function of stimulation frequency, after 40 periods of stimulation.  Colors correspond to the neurons' initial phases.  \label{type2_phase_vs_freq}}
\end{figure}

We simulated 500 Hodgkin-Huxley neurons with initial phases evenly spaced between $0$ and $2 \pi$, corresponding to an initial uniform phase distribution.  The stimulation frequency was varied from 70 Hz to 300 Hz in increments of 5 Hz.  Figure~\ref{type2_phase_vs_freq} shows the final phases after 40 periods of stimulation, after transients have decayed away.  The colors indicate the initial phases of the neurons.  Not all colors are visible for most stimulation frequencies because the final phases of entire subpopulations of neurons are nearly identical, and only one representative initial phase can be seen.  All of the neurons which have nearly the same final phase are part of the same cluster.
 
\begin{figure}[t]
\begin{center}
\leavevmode
\epsfxsize=3.2in
\epsfbox{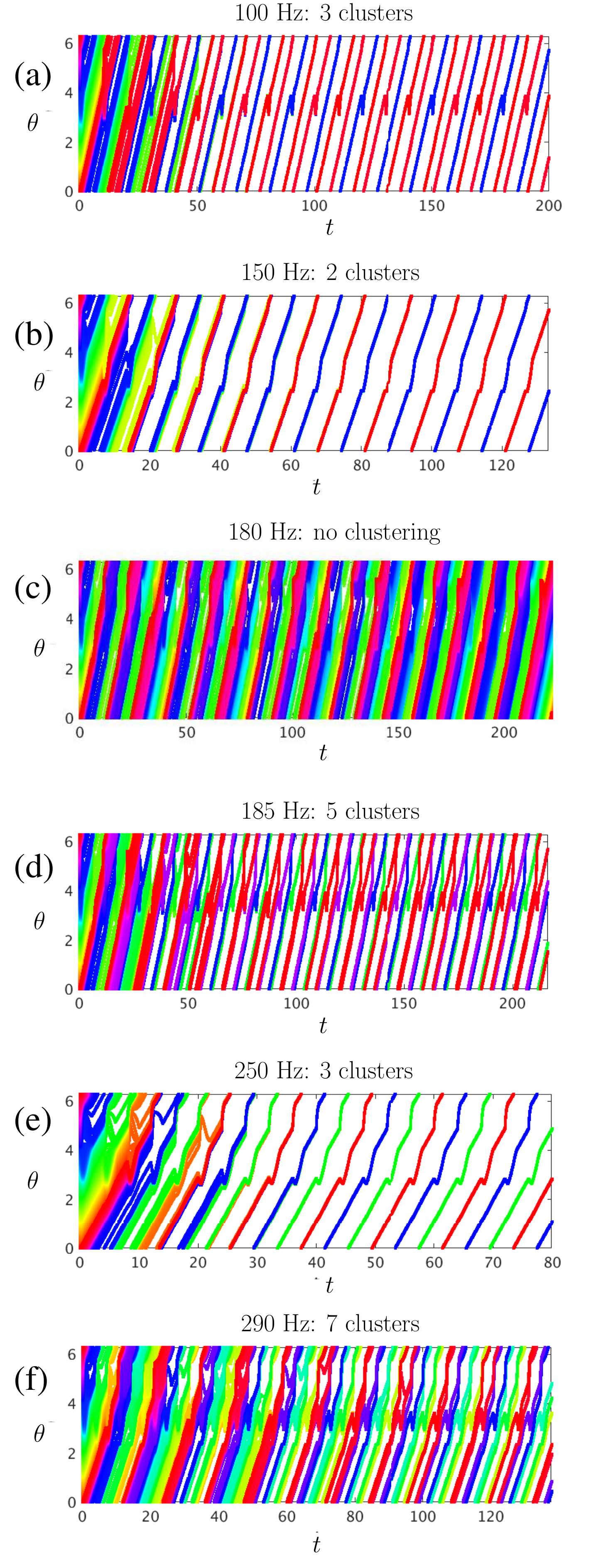}
\end{center}
\caption{Time series showing the phases of Hodgkin-Huxley neurons drawn from an initial uniform distribution for frequencies (a) 100 Hz, (b) 150 Hz, (c) 180 Hz, (d) 185 Hz, (e) 250 Hz, and (f) 290 Hz.  The titles of these panels indicate the number of clusters found after transients have decayed away.  For (c), clusters do not form.  For this and subsequent time series figures, $t$ is measured in ms, and the colors indicate the initial phases of the neurons, with colorbar as in Figure~\ref{type2_phase_vs_freq}. \label{phase_versus_time}}
\end{figure}

Figure~\ref{phase_versus_time} shows the times series of the phases of a population of Hodgkin-Huxley neurons for selected frequencies, and helps us to interpret the results shown in Figure~\ref{type2_phase_vs_freq}.  For example, Figure~\ref{phase_versus_time}(a) shows that for a $100$ Hz stimulus the neurons separate into three clusters, as is also the case for $250$ Hz as shown in Figure~\ref{phase_versus_time}(e).  (Notice that Figure~\ref{type2_phase_vs_freq} shows three possible final phases for each of these frequencies, corresponding to these three clusters.)  Figure~\ref{phase_versus_time}(b) shows that for a $150$ Hz stimulus they separate into two clusters.  For a $180$ Hz stimulus, there is no clustering; see Figure~\ref{phase_versus_time}(c).  By carefully looking at Figure~\ref{phase_versus_time}(d), one sees that for a $185$ Hz stimulus there are five clusters, and from Figure~\ref{phase_versus_time}(f) that for a $295$ Hz stimulus there are seven clusters, as expected from final states shown at these frequencies in Figure~\ref{type2_phase_vs_freq}.  Such clustering behavior and non-clustering (chaotic) behavior has been seen in other studies, such as~\cite{wils11} and~\cite{wils15cluster}.

\begin{figure}[tb]
\begin{center}
\leavevmode
\epsfxsize=3.2in
\epsfbox{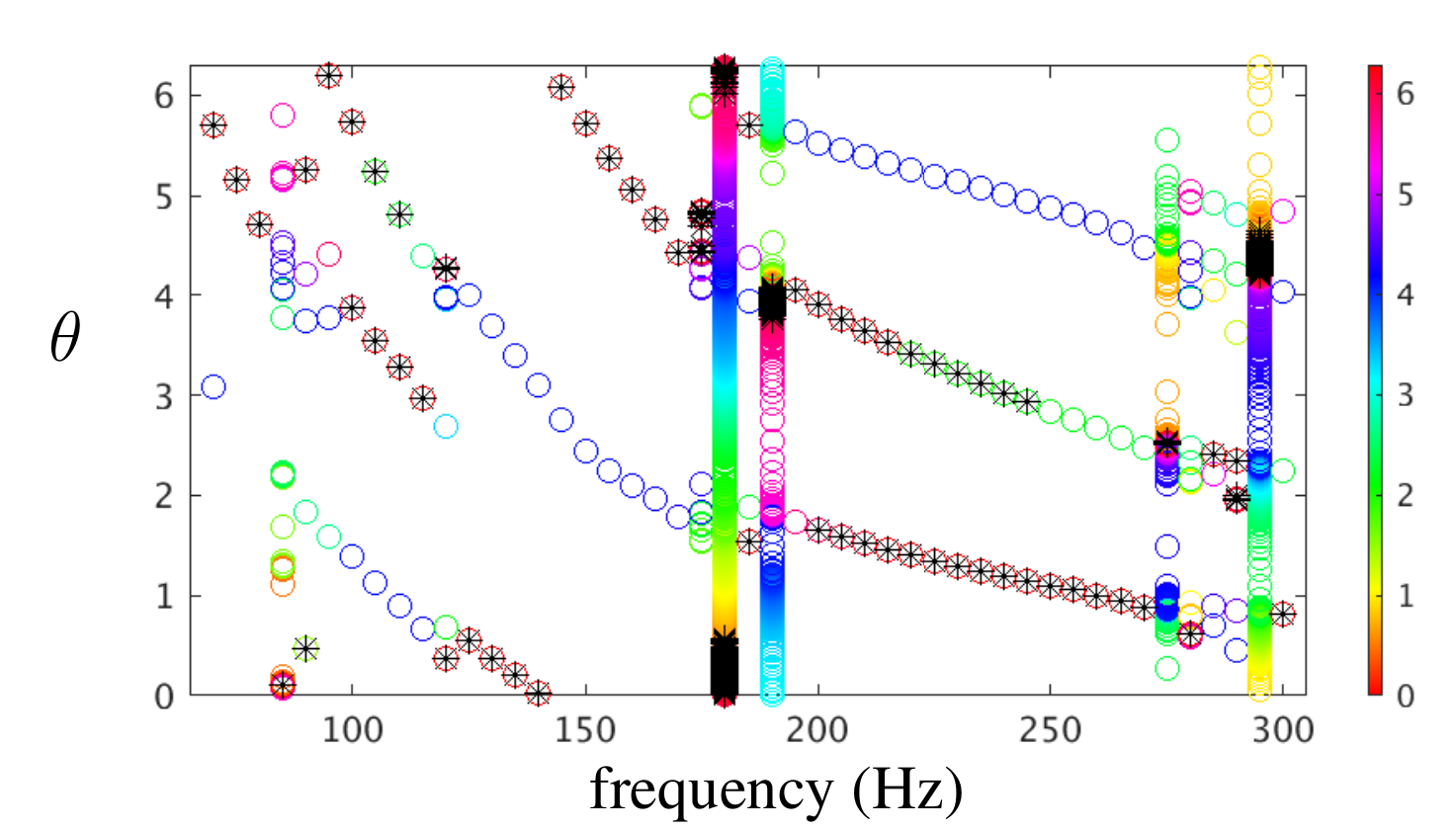}
\end{center}
\caption{As a function of stimulation frequency, the final phases $\theta$ of Hodgkin-Huxley neurons drawn from an initial von Mises distribution after 40 periods of stimulation are shown as black $*$'s, overlaid on the final phases of Hodgkin-Huxley neurons drawn from an initial uniform distribution (as was shown in Figure~\ref{type2_phase_vs_freq}). \label{type2_phase_vs_freq_von_mises}}
\end{figure}

Inspired by neural synchrony in Parkison's patients, we also considered an initial partially synchronized neural population, with phases distributed according to a von Mises distribution~\cite{best79} centered at $\theta = 0$:
\begin{equation}
\rho_0(\theta) = \frac{e^{\kappa \cos \theta}}{2 \pi I_0 (\kappa)},
\end{equation}
where $I_0(\kappa)$ is the modified Bessel function of order 0.  This distribution is similar to a Gaussian distribution, but on a circle.  
We simulated 500 Hodgkin-Huxley neurons with initial phases distributed according to the von Mises distribution with $\kappa = 50$.  As for Figure~\ref{type2_phase_vs_freq}, the stimulation frequency was varied from 70 Hz to 300 Hz in increments of 5 Hz.  Figure~\ref{type2_phase_vs_freq_von_mises} shows the final phases after 40 periods of stimulation, after transients have decayed away.  We see that the final phases of the neurons from the initial von Mises distribution lie on a subset of the final phases of the neurons from the intial uniform distribution.  For example, when the stimulation frequency is 100 Hz, the neurons from the initial von Mises distribution are concentrated in two of the three clusters which exist for the initial uniform distribution.

\clearpage

\begin{figure}[tb!]
\begin{center}
\leavevmode
\epsfxsize=2.2in
\epsfbox{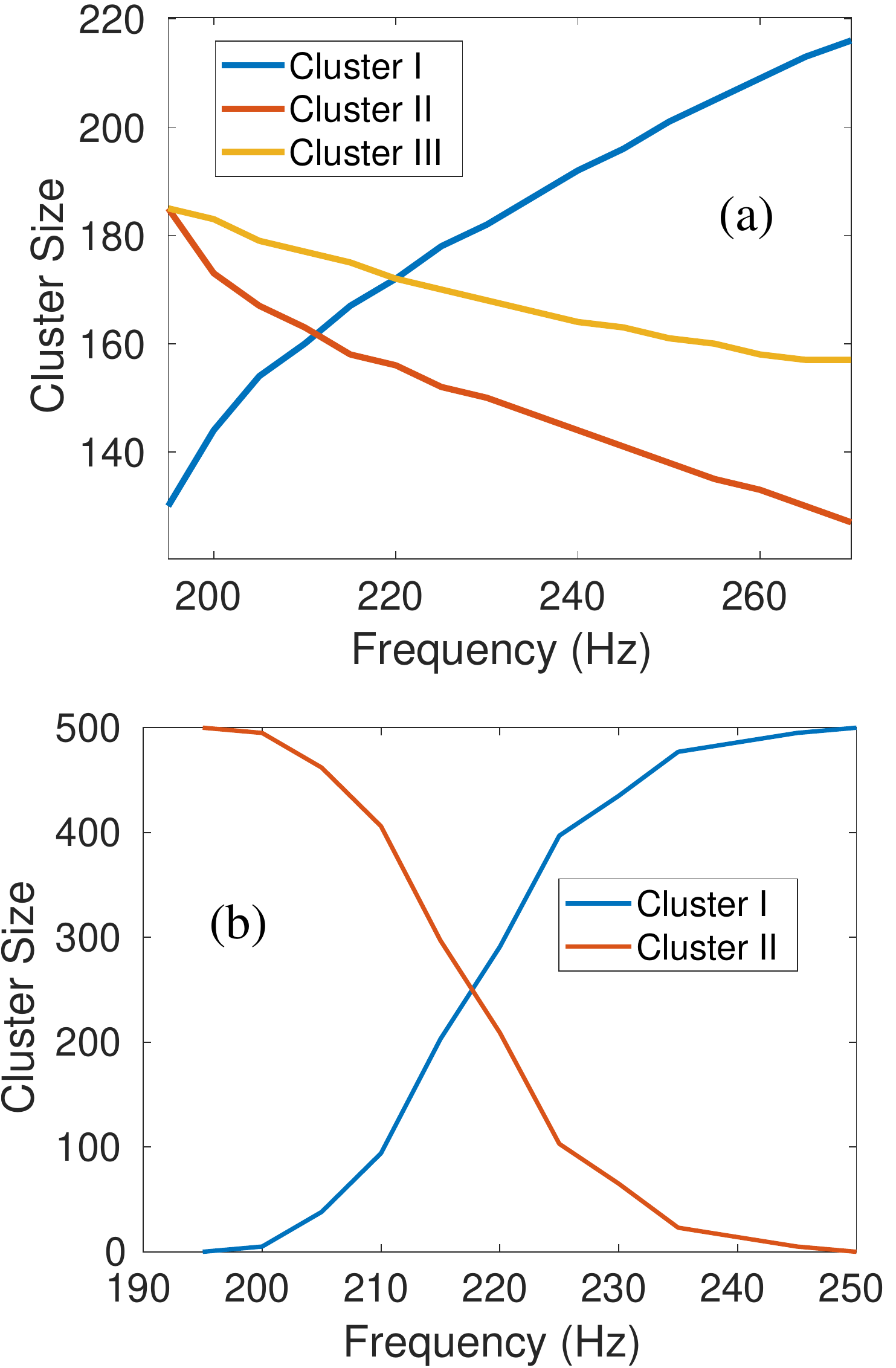}
\end{center}
\caption{The number of Hodgkin-Huxley neurons in different clusters for a population size of 500, with initial (a) uniform and (b) von Mises distributions.  \label{cluster_sizes}}
\end{figure}

We also designed an algorithm to detect the size of clusters in a population. The algorithm groups the phases of neurons in a population at each timestep into clusters by sorting the phases in ascending order and checking if the $i^{\rm th}$ phase is within $\epsilon$ of the $(i + 1)$-th phase for an appropriate small value of $\epsilon$. If so, the size of the current cluster is increased by one. If not, the algorithm creates a new cluster. The process is repeated until all neurons have been grouped into clusters.  Figure~\ref{cluster_sizes} shows the number of neurons in the different clusters over a range of frequencies for the intial uniform distribution (for which three clusters are populated) and von Mises distribution (for which only two clusters are populated).  As we will see in Section~\ref{section:identical_stimuli}, this figure can be explained in terms of the basins of attraction of fixed points of iterates of a one-dimensional map defined on the circle.  The initial phase of a given neuron will determine which cluster it ends up in.

\begin{figure}[tb]
\begin{center}
\leavevmode
\epsfxsize=3.2in
\epsfbox{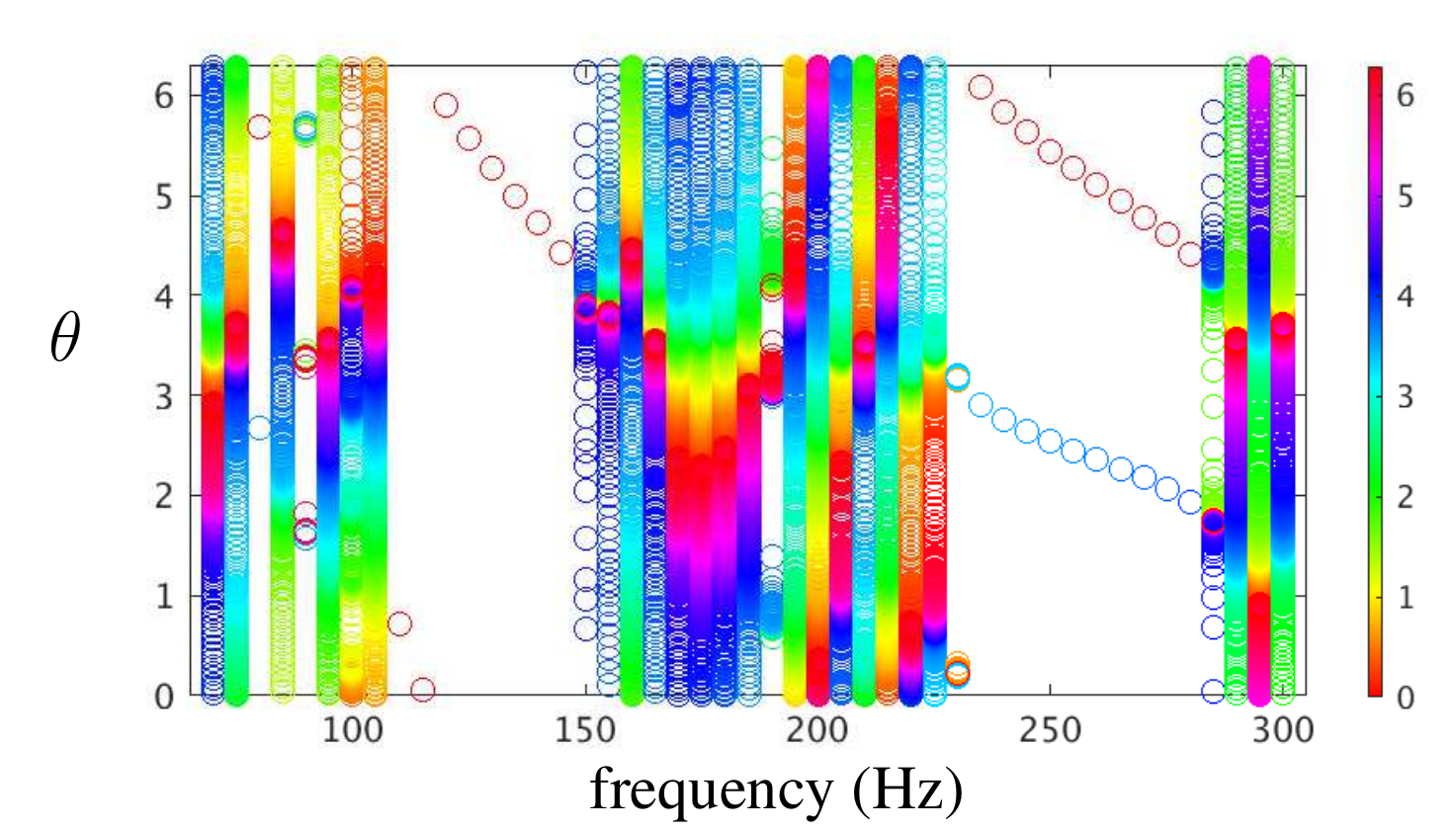}
\end{center}
\caption{The final phases $\theta$ of thalamic neurons drawn from an initial uniform distribution as a function of stimulation frequency, after 40 periods of stimulation. Colors correspond to the neurons' initial phases.  \label{type1_phase_vs_freq}}
\end{figure}

\begin{figure}[tb!]
\begin{center}
\leavevmode
\epsfxsize=3.5in
\epsfbox{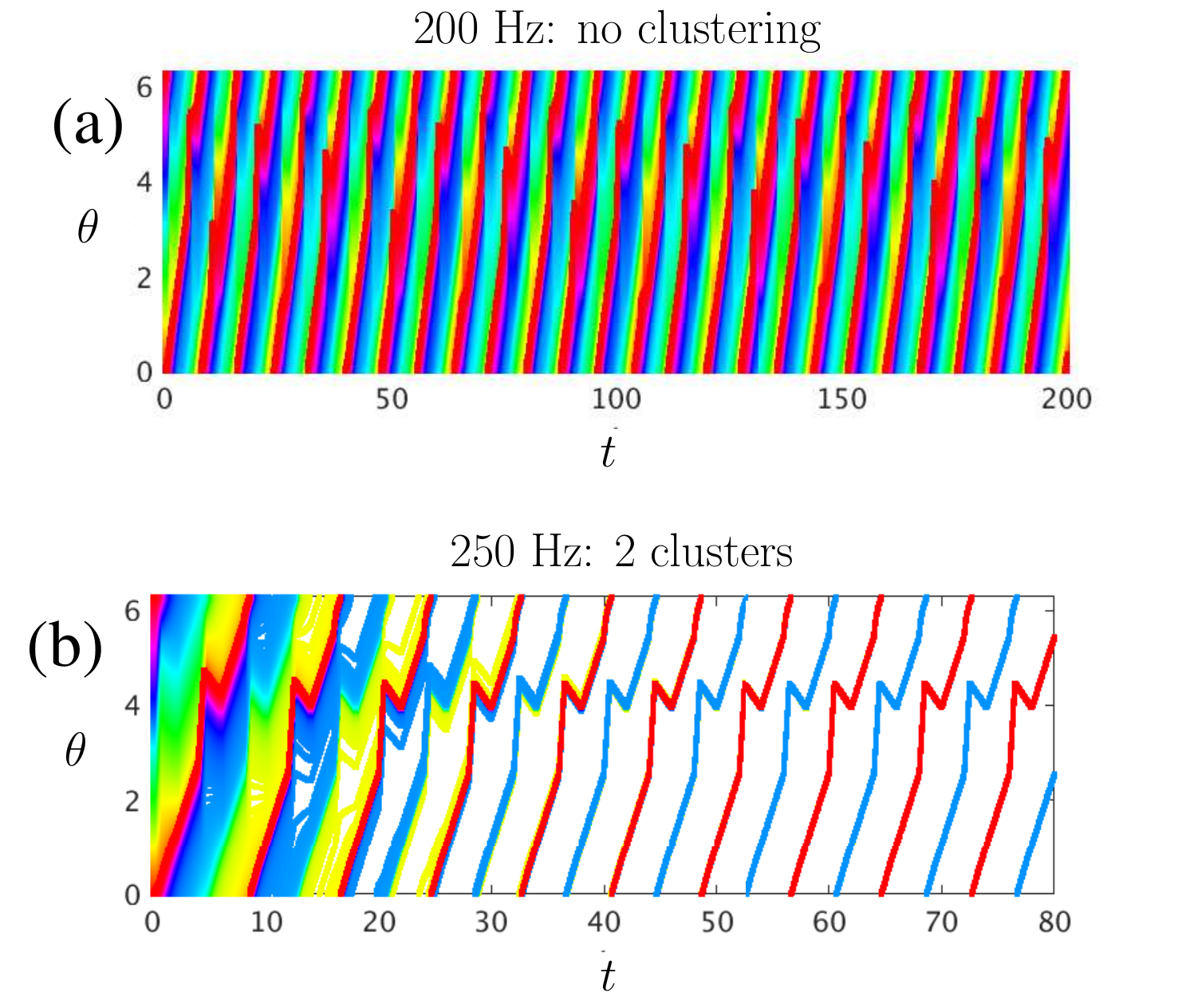}
\end{center}
\caption{Time series showing the phases of thalamic neurons drawn from an initial uniform distribution for frequencies (a) 200 Hz, and (b) 250 Hz. For (a), clusters do not form; for (b), there are two clusters after transients decay away. \label{phase_versus_time_type1}}
\end{figure}

We also considered populations of thalamic (Type I) neurons with the same stimuli~(\ref{eq:Pulsatile_Stimulus}) with $u_{max}$ corresponding to a current density of $20$ $\mu A/cm^2$, $p = 0.5$ ms, and $\lambda = 3$.  
We simulated 500 thalamic neurons with initial phases evenly spaced between $0$ and $2 \pi$, corresponding to a uniform distribution.  The stimulation frequency was varied from 70 Hz to 300 Hz in increments of 5 Hz.  Figure~\ref{type1_phase_vs_freq} shows the final phases after 40 periods of stimulation, after transients have decayed.  Figure~\ref{phase_versus_time_type1} shows the time series of the phases of a population of such neurons for selected frequencies.  Here we again see clustering for some frequencies (such as $250$ Hz), and non-clustering behavior for other frequencies (such as $200$ Hz).

In the next section, we derive and investigate one-dimensional maps which can be used to understand the types of clusters which occur in these simulations, along with their stability properties and their basins of attraction. 

\section{Analysis of Clusters due to Identical Pulses} \label{section:identical_stimuli}

In this section, we show how the clustering behavior found in the simulations from Section~\ref{section:simulations} can be understood in terms of appropriate compositions of one-dimensional maps on the circle.  

\begin{figure}[h!]
\begin{center}
\leavevmode
\epsfxsize=2.2in
\epsfbox{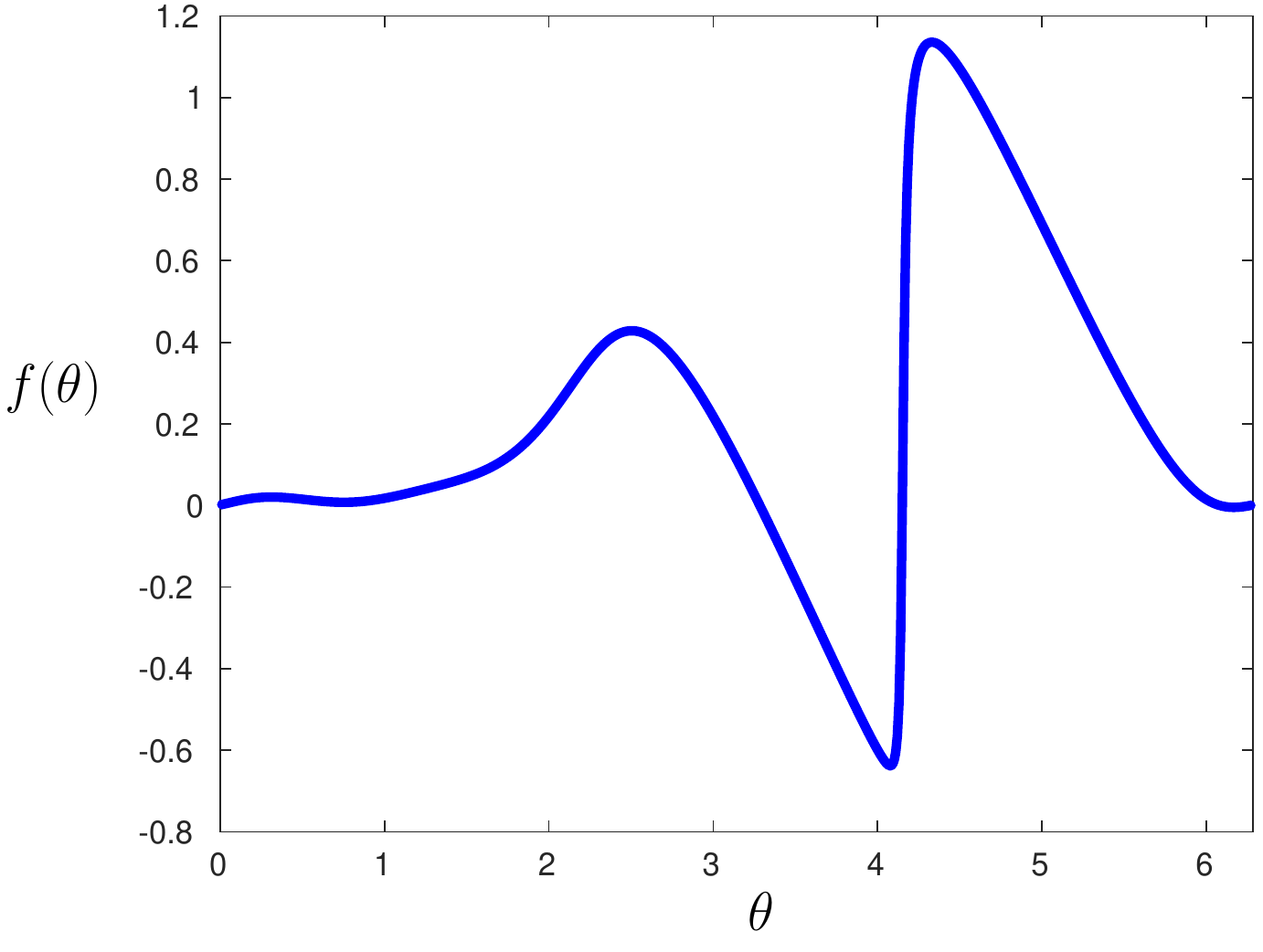}
\end{center}
\caption{Response function $f(\theta)$ which characterizes the phase response of Hodgkin-Huxley neurons to the stimulus, for $u_{max}$ corresponding to a current density of $20 \mu A/cm^2$, $p = 0.5$ ms, and $\lambda=3$.
\label{f_type2}}
\end{figure}

We consider a system of neural oscillators subjected to a $\tau$-periodic sequence of pulses as shown in Figure~\ref{one_kick}, and described by the dynamics~\cite{wils15cluster}
\begin{equation}
\dot{\theta}_i = \omega + f(\theta_i) \delta({\rm mod}(t,\tau)), \qquad i = 1,\cdots,N.
\end{equation}
Here the response function $f(\theta)$ describes the change in phase due to a single pulse (including the positive current for time $p$, and the negative current for time $\lambda p$).  If the pulse was a delta function with unit area, $f(\theta)$ would be equal to the infinitesimal PRC $Z(\theta)$; for more general pulses, it can be calculated using a direct method in which a pulse is applied at a known phase, and the change in phase is deduced from the change in timing of the next action potential~\cite{neto12}.  We will think of the change in phase due to the pulse as occurring instantaneously, even though the pulse will typically have a finite duration; this will be a good approximation for pulses of short duration.  Figure~\ref{f_type2} shows $f(\theta)$ for the Hodgkin-Huxley neurons considered in this paper for pulses as shown in Figure~\ref{one_kick} with $u_{max}$ corresponding to a current density of $20 \mu A/cm^2$, $p = 0.5$ ms, and $\lambda = 3$.

To understand the clustering behavior, it will be useful to consider the map which takes the phase of a neuron to the phase exactly one forcing cycle later, cf.~\cite{wils15cluster}.  To find this map, suppose that we start with $\theta(0^+) = 0$, immediately after the start of a pulse, where we assume that we have already accounted for the effect of the pulse according to the function $f(\theta)$. The next pulse comes at time $\tau$.  Up until time $\tau$, the phase evolves according to $\dot{\theta} = \omega$; therefore,
\begin{equation}
\theta(\tau^-) = \theta_0 + \omega \tau.
\end{equation}
Treating the change in phase due to the next pulse as occurring instantaneously, we have
\begin{equation}
\theta(\tau^+) = \theta_0 + \omega \tau + f(\theta_0 + \omega \tau).
\end{equation}
The system then evolves for a time $\tau$ without stimulus, giving
\begin{equation}
\theta(2 \tau^-) = \theta_0 + 2 \omega \tau + f(\theta_0 + \omega \tau);
\end{equation}
the next pulse at time $2 \tau$ gives
\begin{equation}
\theta(2 \tau^+) = \theta_0 + 2 \omega \tau + f(\theta_0 + \omega \tau) + f(\theta + 2 \omega \tau + f(\theta_0 + \omega \tau)),
\end{equation}
and so on.  It is useful to let~\cite{wils15cluster}
\begin{equation}
g(s) = s + \omega \tau + f(s + \omega \tau),
\end{equation}
which gives
\begin{equation}
\theta(n \tau^+) = g^{(n)} (\theta_0), 
\end{equation}
where $g^{(n)}$ denotes the composition of $g$ with itself $n$ times, and $\theta_0$ is the initial state of the neuron.  

We look for fixed points of $g^{(n)}$, that is, solutions to $\theta^* = g^{(n)}(\theta^*)$; for such solutions, the phase has the same value after $n$ pulses as where it started.   We are particularly interested in fixed points of $g^{(n)}$ which are not fixed points of $g^{(m)}$ for any positive integer $m$ satisfying $m < n$; then there will be $n$ fixed points of $g^{(n)}$ that correspond to points on a period-$n$ orbit of $g$.  If 
\begin{equation}
\left| \left. \frac{d}{d \theta} \right|_{\theta = \theta^*} (g^{(n)} (\theta)) \right|<1,
\end{equation}
then the fixed point $\theta^*$ of $g^{(n)}$ is stable, as is the corresponding period-$n$ orbit of $g$.  Neurons which start with initial phases within the basin of attraction of a given fixed point of $g^{(n)}$ will asymptotically approach that fixed point under iterations of $g^{(n)}$.  The $n$ different fixed points will each have a basin of attraction, so a uniform intial distribution of neurons will form $n$ clusters, one for each of these fixed points of $g^{(n)}$, cf.~\cite{wils15cluster}.  

\begin{figure}[t!]
\begin{center}
\leavevmode
\epsfxsize=2.2in
\epsfbox{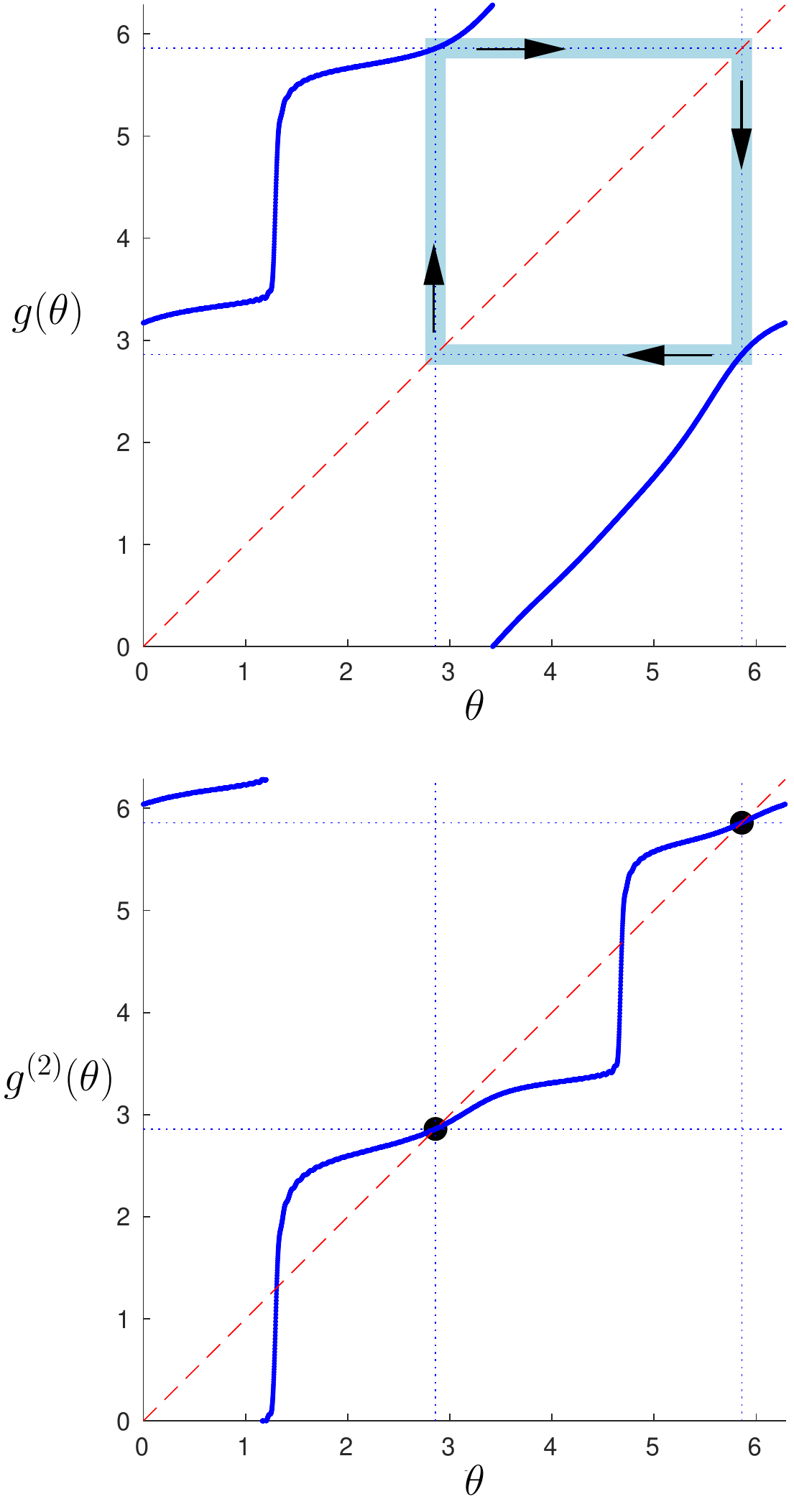}
\end{center}
\caption{Maps $g(\theta)$ and $g^{(2)}(\theta)$ for Hodgkin-Huxley neuron with stimulation frequency $150$ Hz.  Intersections with the diagonal dashed line indicate fixed points of the respective map.  The dotted lines show $\theta$ values for the stable fixed points of the $g^{(2)}$ map.  
\label{f150}}
\end{figure}

We now illustrate how these maps can be used to understand the specific clustering behavior shown in Section~\ref{section:simulations}.
As a first example, suppose that a population of Hodgkin-Huxley neurons is stimulated with frequency 150 Hz, corresponding to $\tau = 6.67$ ms.  Figure~\ref{f150} shows $g(\theta)$ and $g^{(2)}(\theta)$.  Fixed points of these maps correspond to intersections with the diagonal.   We see that there are two stable fixed points for $g^{(2)}(\theta)$, at $\theta = 2.86$ and $\theta = 5.86$ (these fixed points are stable because the slope at the intersection is between $-1$ and $+1$).  There are also two unstable fixed points for $g^{(2)}$ at $\theta = 1.305$ and $\theta = 4.685$, where the slope at the intersection is greater than 1.  There are no fixed points for $g(\theta)$, but a cobweb analysis verifies that there is a period-2 orbit
\[
\theta = 2.86 \rightarrow 5.86 \rightarrow 2.86 \rightarrow \cdots.
\]
These fixed points of $g^{(2)}$ correspond a stable 2-cluster state for a population of oscillators, as shown in Figure~\ref{phase_versus_time}(b).  We note that we can deduce the basin of attraction for the different stable fixed points of $g^{(2)}$; for example, the basin of attraction for the stable fixed point at $\theta = 2.86$ is the range $1.305 < \theta_0 < 4.685$, that is, between the two unstable fixed points.

\begin{figure}[b!]
\begin{center}
\leavevmode
\epsfxsize=2.2in
\epsfbox{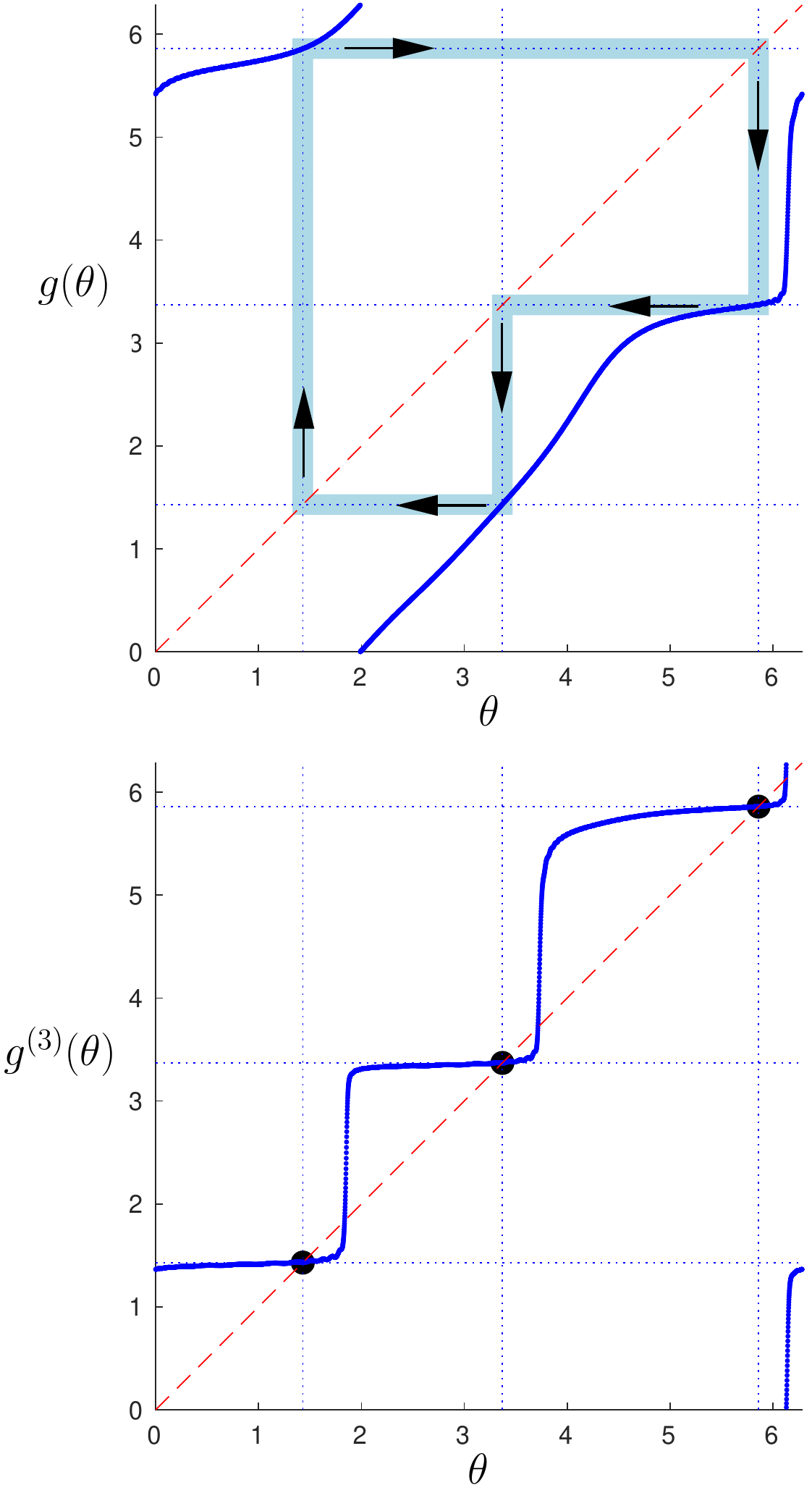}
\end{center}
\caption{Maps $g(\theta)$ and $g^{(3)}(\theta)$ for Hodgkin-Huxley neuron with stimulation frequency $100$ Hz.
\label{f100}}
\end{figure}

As another example, suppose that a population of Hodgkin-Huxley neurons is stimulated with frequency 100 Hz, corresponding to $\tau = 10$ ms.  Figure~\ref{f100} shows $g(\theta)$ and $g^{(3)}(\theta)$.  We see that there are three stable fixed points for $g^{(3)}(\theta)$, at $\theta = 1.43$, $\theta = 3.37$, and $\theta = 5.86$ (these fixed points are stable because the slope at the intersection has slope between $-1$ and $+1$).  There are no fixed points for $g(\theta)$, but a cobweb analysis verifies that there is a period-3 orbit
\[
\theta = 1.43 \rightarrow 5.86 \rightarrow 3.37 \rightarrow 1.43 \rightarrow \cdots.
\]
These fixed points of $g^{(3)}$ correspond a stable 3-cluster state for a population of oscillators, as shown in Figure~\ref{phase_versus_time}(a).


\begin{figure}[t!]
\begin{center}
\leavevmode
\epsfxsize=2.2in
\epsfbox{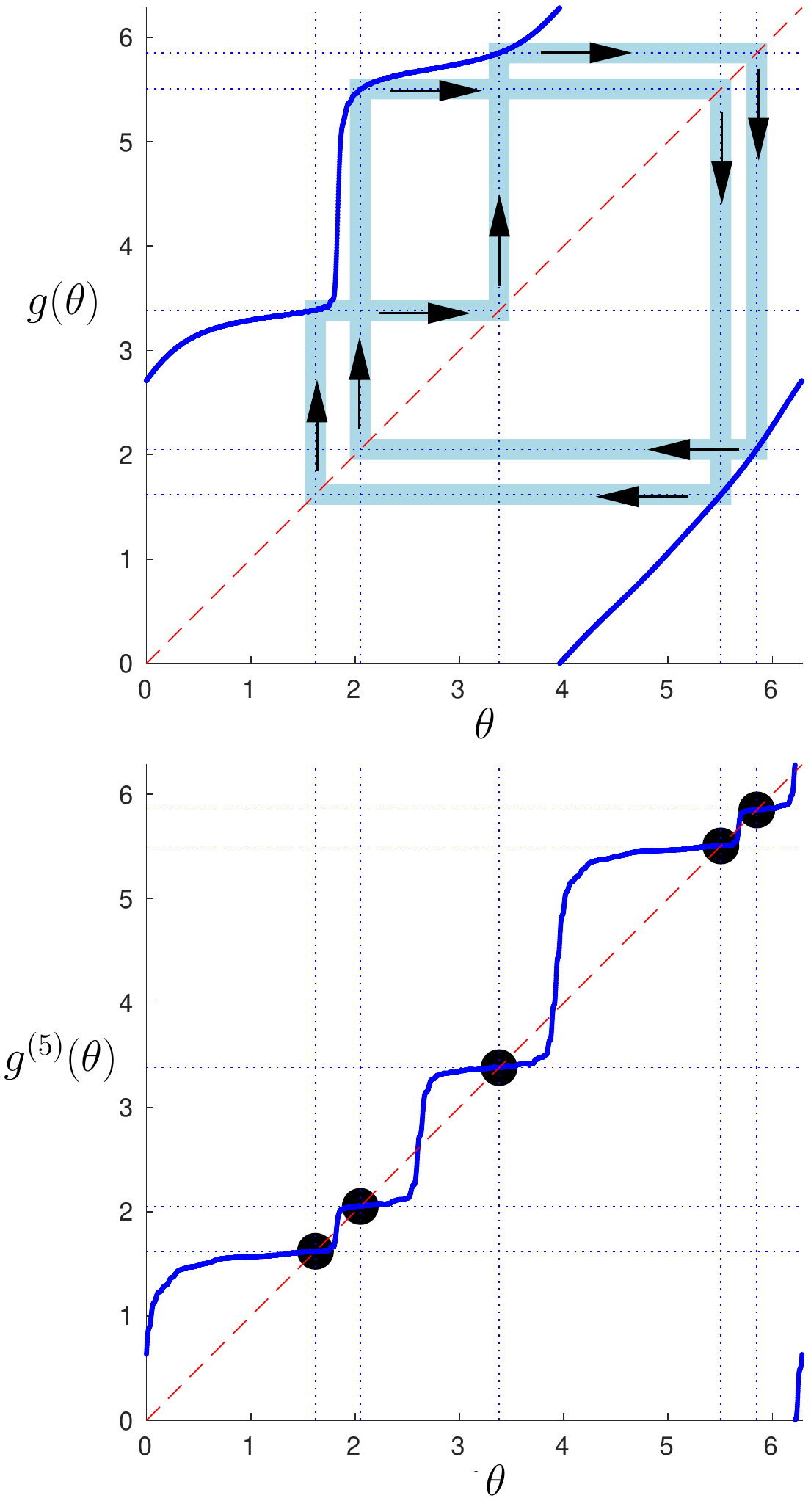}
\end{center}
\caption{Maps $g(\theta)$ and $g^{(5)}(\theta)$ for Hodgkin-Huxley neuron with stimulation frequency $185$ Hz.
\label{f185}}
\end{figure}

We see that this map can also capture $n$-clusters for larger values of $n$.  For example, for frequency 185 Hz, the $g$ and $g^{(5)}$ maps shown in Figure~\ref{f185} confirm that there is a stable period-5 orbit
\[
\theta = 1.62 \rightarrow 3.38 \rightarrow 5.85 \rightarrow 2.05 \rightarrow 5.51 \rightarrow 1.62 \rightarrow \cdots,
\]
corresponding to the stable 5-cluster state shown in Figure~\ref{phase_versus_time}(d).

\begin{figure}[h!]
\begin{center}
\leavevmode
\epsfxsize=2.2in
\epsfbox{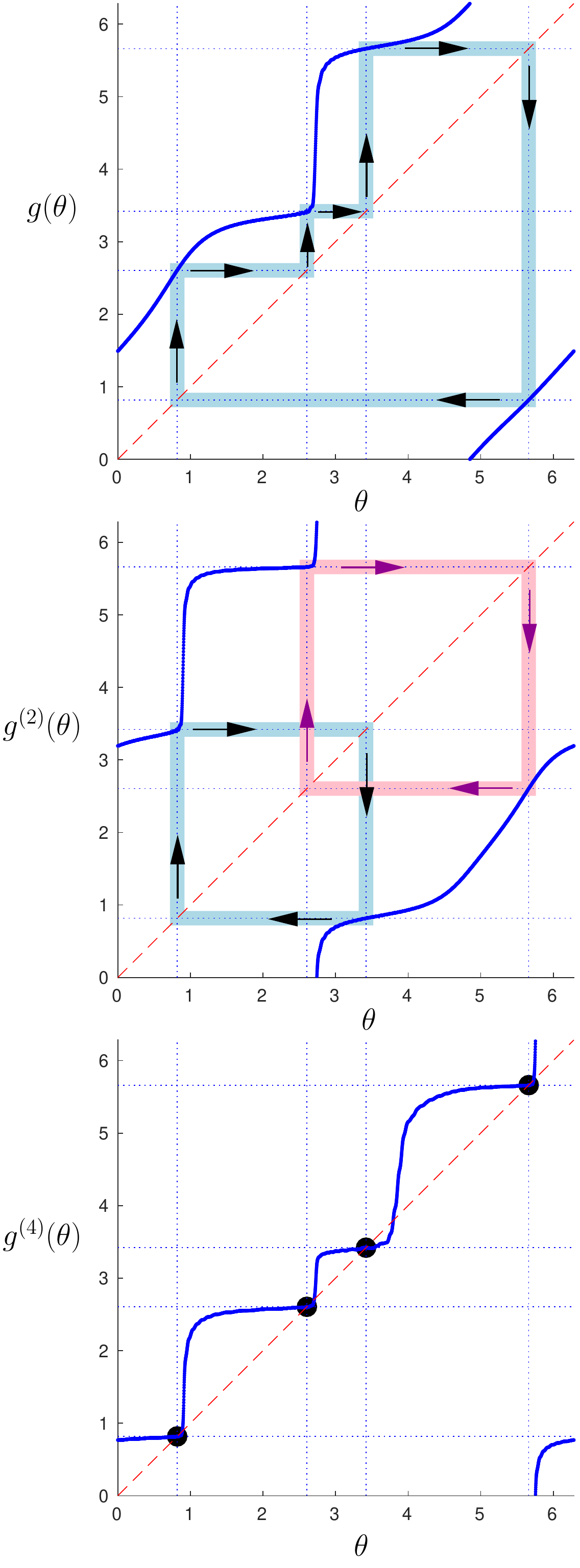}
\end{center}
\caption{For the Hodgkin-Huxley neurons with stimulation frequency $300$ Hz, there is a stable period-4 orbit for $g$, which corresponds to two stable period-2 orbits for $g^{(2)}$, which in turn correspond to four stable fixed points for $g^{(4)}$. 
\label{f300}}
\end{figure}

Finally, if we apply these identical stimuli at $300$ Hz,
we obtain a total of four stable fixed points for $g^{(4)}$, corresponding to a stable period-4 orbit for $g(\theta)$ 
\[
\theta = 0.82 \rightarrow 2.61 \rightarrow 3.32 \rightarrow 5.69 \rightarrow 0.82 \rightarrow \cdots\
\]
which is equivalent to two stable period-2 orbits for $g^{(2)}$
\[
\theta = 0.82 \rightarrow 3.32 \rightarrow 0.82 \rightarrow \cdots, 
\]
\[
\theta = 2.61 \rightarrow 5.69 \rightarrow 2.61 \rightarrow \cdots, 
\]
and four stable fixed points for $g^{(4)}$
\[
\theta = 0.82, \qquad \theta = 2.61, \qquad \theta = 3.32, \qquad \theta = 5.69;
\]
see Figure~\ref{f300}.  We show in the next section that it is possible to obtain similar dynamics with stimuli consisting of pulses with alternating properties.

\begin{figure}[tb]
\begin{center}
\leavevmode
\epsfxsize=2.2in
\epsfbox{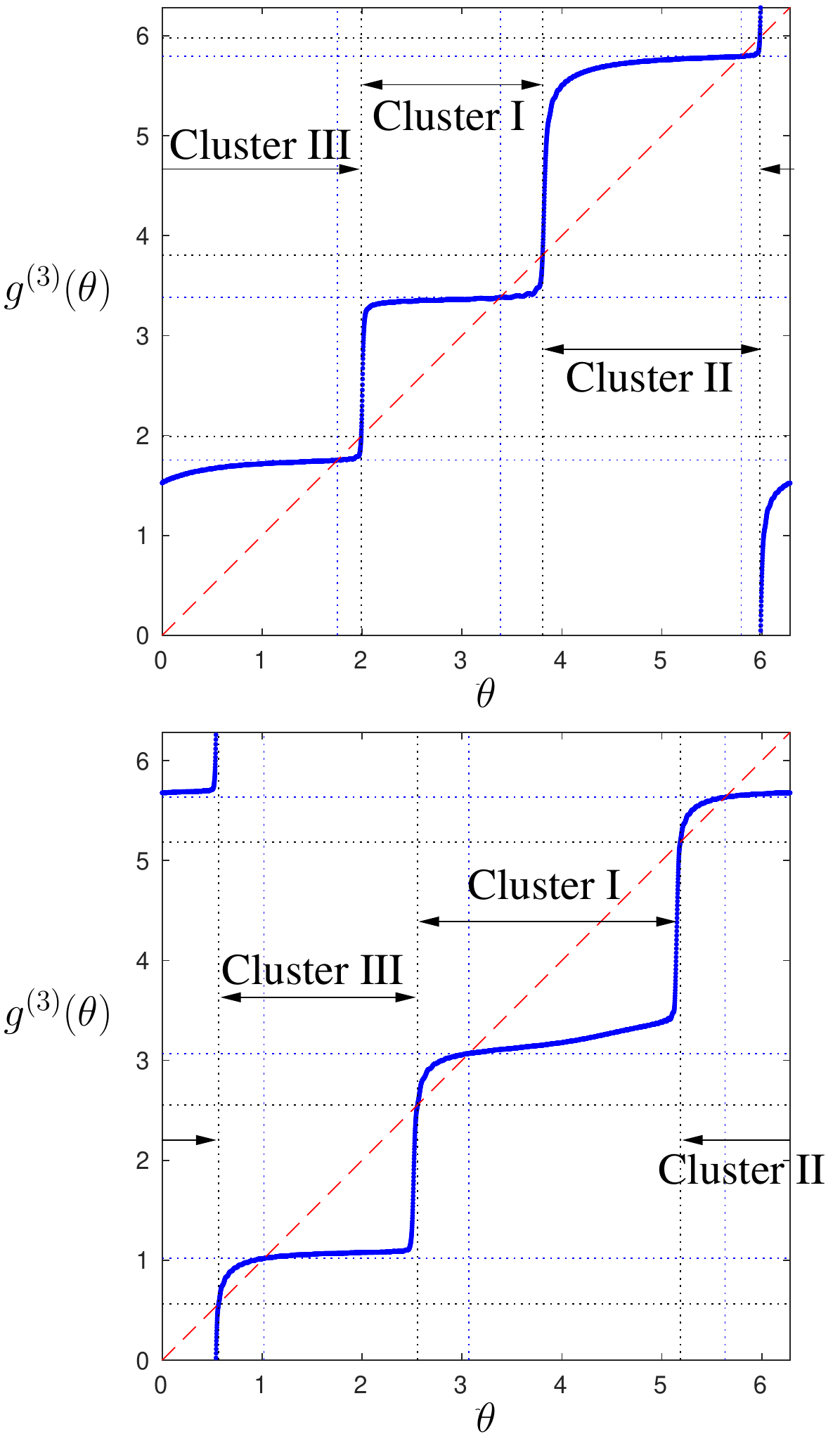}
\end{center}
\caption{Basins of attraction for the different clusters for (a) $200$ Hz and (b) $260$ Hz stimuli.
\label{cluster_sizes_map}}
\end{figure}

We can understand the cluster sizes shown in Figure~\ref{cluster_sizes}(a) by looking at the basins of attraction of the different stable fixed points, as indicated in Figure~\ref{cluster_sizes_map} for $200$ Hz and $260$ Hz stimuli.  The basin boundaries are at the phases of the appropriate unstable fixed points.  When the initial phase distribution is uniform, the number of neurons which end up in each cluster is proportional to the size of the corresponding basin of attraction.  For example, if there are 500 uniformly distributed neurons, this predicts that there will be 144, 173, and 183 neurons in Clusters I, II, and III, respectively, for a $200$ Hz stimluus, and 209, 133, and 159 neurons in Clusters I, II, and III, respectively, for a $260$ Hz stimulus.  This is consistent with the results shown in Figure~\ref{cluster_sizes}(a).  The number of neurons in each cluster for Figure~\ref{cluster_sizes}(b) would be determined by the number of neurons which are initially in the respective basin of attraction, as determined by the initial phase distribution; here, there were no neurons with initial phases that end up in Cluster III.

\begin{figure}[tb!]
\begin{center}
\leavevmode
\epsfxsize=2.2in
\epsfbox{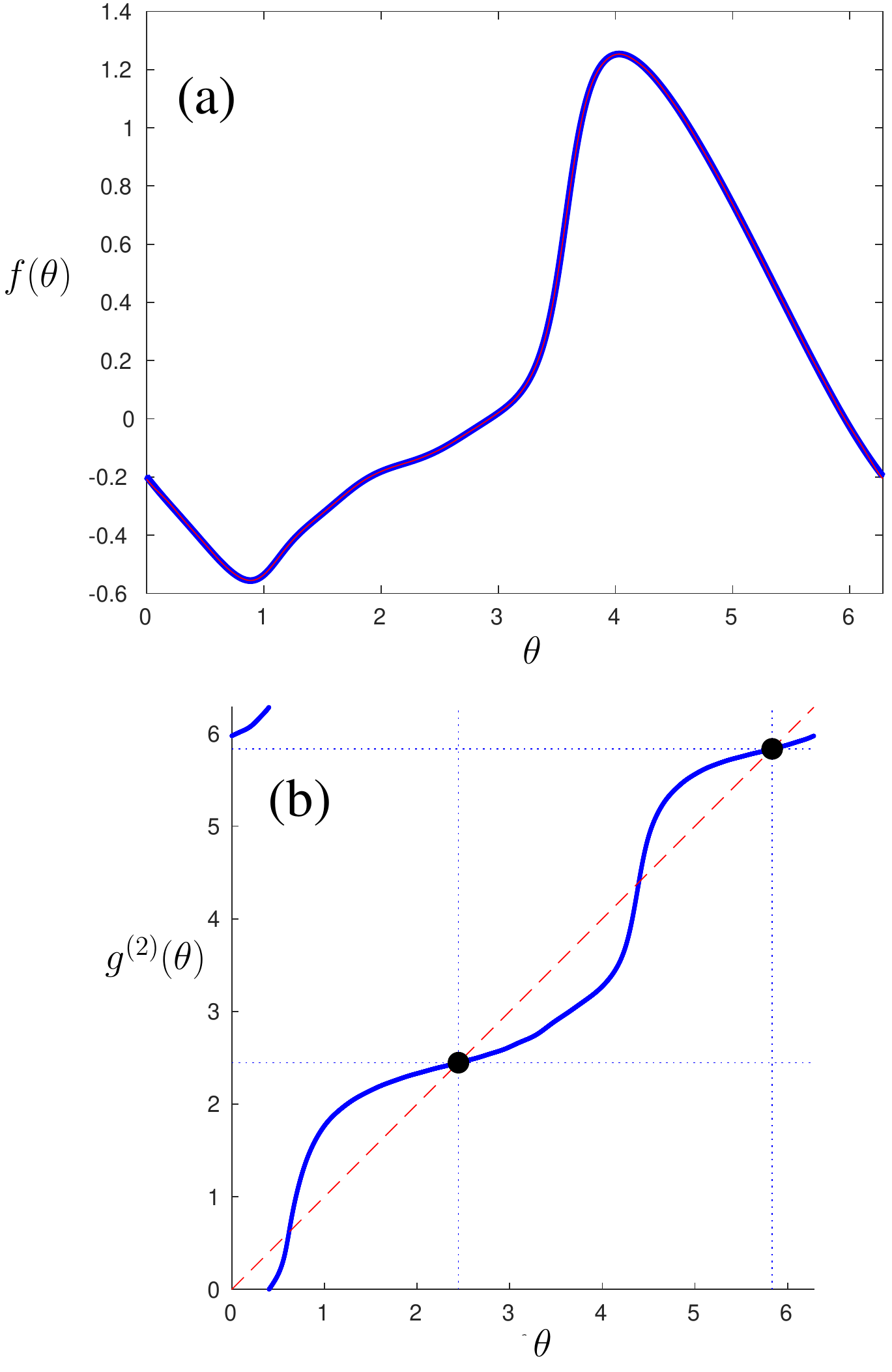}
\end{center}
\caption{(a) Response function $f(\theta)$ which characterizes the phase response of thalamic neurons to a pulse with $u_{max}$ corresponding to a current of $20 \mu A/cm^2$, $p = 0.5$ ms, and $\lambda=3$.  (b) Map $g^{(2)}(\theta)$ for the thalamic neuron with stimulation frequency $250$ Hz, showing two stable fixed points which correspond to a 2-cluster state.
\label{type1}}
\end{figure}

The same analysis techniques can also be used to understand the dynamics of thalamic neurons subjected to periodic pulses.  Figure~\ref{type1}(a) shows the response function $f(\theta)$ for thalamic neurons with the stimulus given by (\ref{eq:Pulsatile_Stimulus}) with $u_{max}$ corresponding to a current density of $20 \mu A/cm^2$, $p = 0.5$ ms, and $\lambda = 3$; Figure~\ref{type1}(b) shows that there is a stable 2-cluster state for a stimulation frequency of $250$ Hz, as expected from Figure~\ref{type1_phase_vs_freq}.

\section{Analysis of Clusters due to Pulses with Alternating Properties} \label{section:alternating_stimuli}

\begin{figure}[h!]
\begin{center}
\leavevmode
\epsfxsize=3.2in
\epsfbox{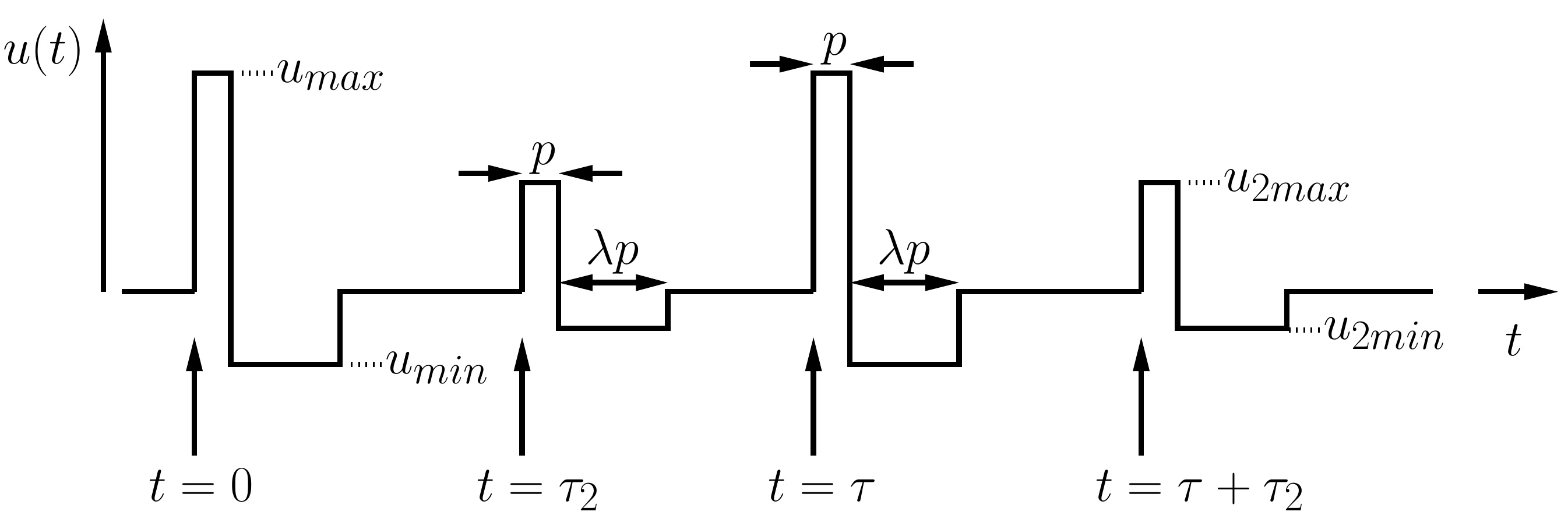}
\end{center}
\caption{Sequence of alternating pulses.
\label{two_kicks}}
\end{figure}

\begin{figure}[h!]
\begin{center}
\leavevmode
\epsfxsize=3.2in
\epsfbox{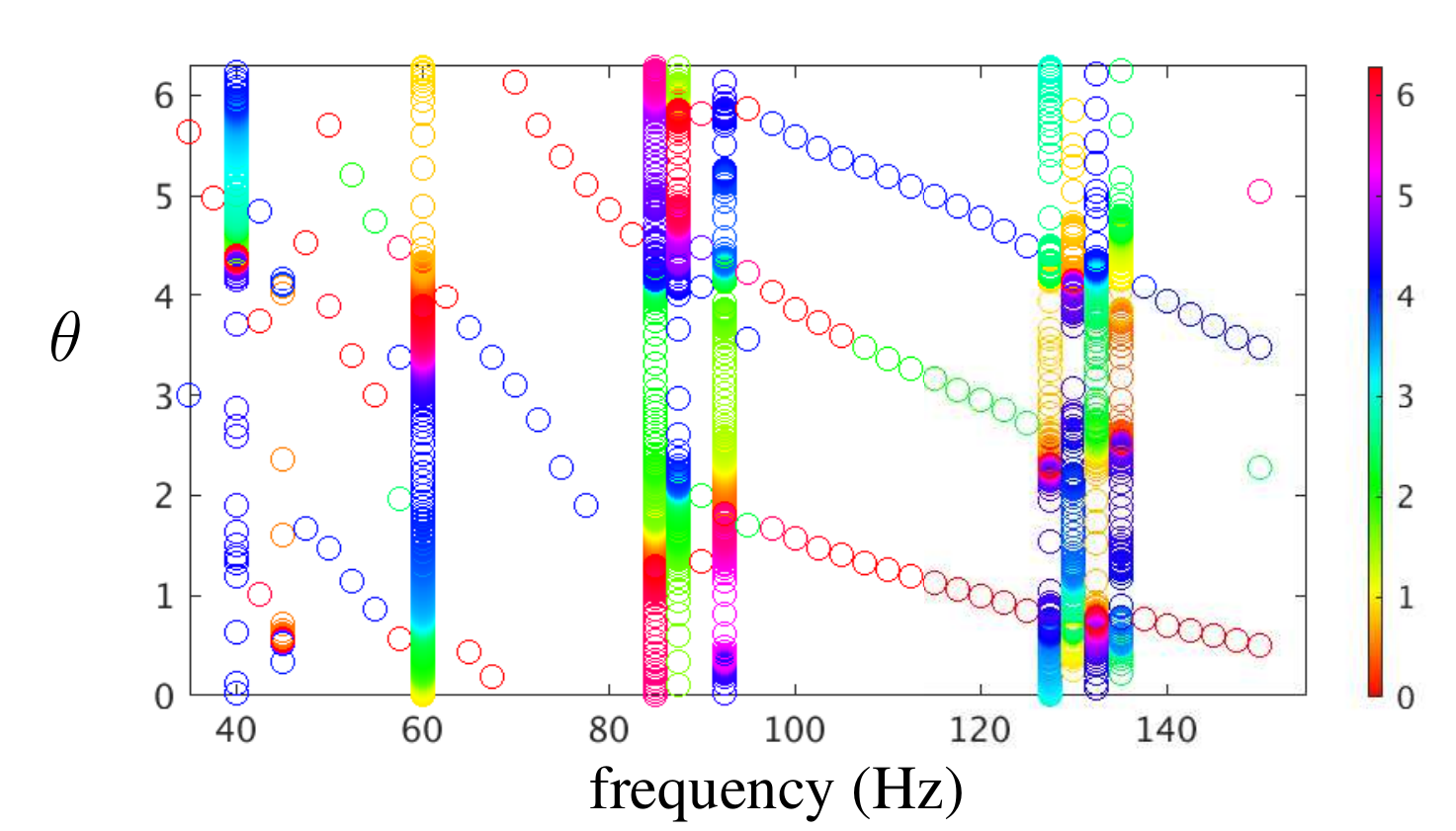}
\end{center}
\caption{The final phases $\theta$ of Hodgkin-Huxley neurons drawn from an initial uniform distribution as a function of stimulation frequency, after 80 periods of pulses with alternating properties (to allow transients to decay), as described in the text. Colors correspond to the neurons' initial phases.  \label{type2_phase_vs_freq_two_kicks}}
\end{figure}

In this section, we consider more general stimuli, specifically pulses with alternating properties, as shown in Figure~\ref{two_kicks}.  Here, the pulses from before, that is with $u_{max}$ corresponding to a current density of $20 \mu A/cm^2$, $p = 0.5$ ms, and $\lambda = 3$, will be assumed to occur at times $0, \tau, 2 \tau, \cdots$.  But now additional pulses with $u_{2max}$ corresponding to a current density of $10 \mu A/cm^2$, $\lambda=3$, $u_{2min} = -u_{2max}/\lambda$ and $p = 0.5$ ms, will be assumed to occur at times $\tau_2$, $\tau + \tau_2$, $2 \tau + \tau_2, \cdots$.  Figure~\ref{type2_phase_vs_freq_two_kicks} shows that the clustering behavior for such alternating pulses with $\tau_2 = \tau/2$ strongly resembles the clustering behavior found at twice the frequency for identical pulses, as shown in Figure~\ref{type2_phase_vs_freq}, although there are differences.  The analysis in this section shows how the methods from Section~\ref{section:identical_stimuli} can be adapted to understand clustering behavior for such alternating pulses.
 
\begin{figure}[ht]
\begin{center}
\leavevmode
\epsfxsize=2.2in
\epsfbox{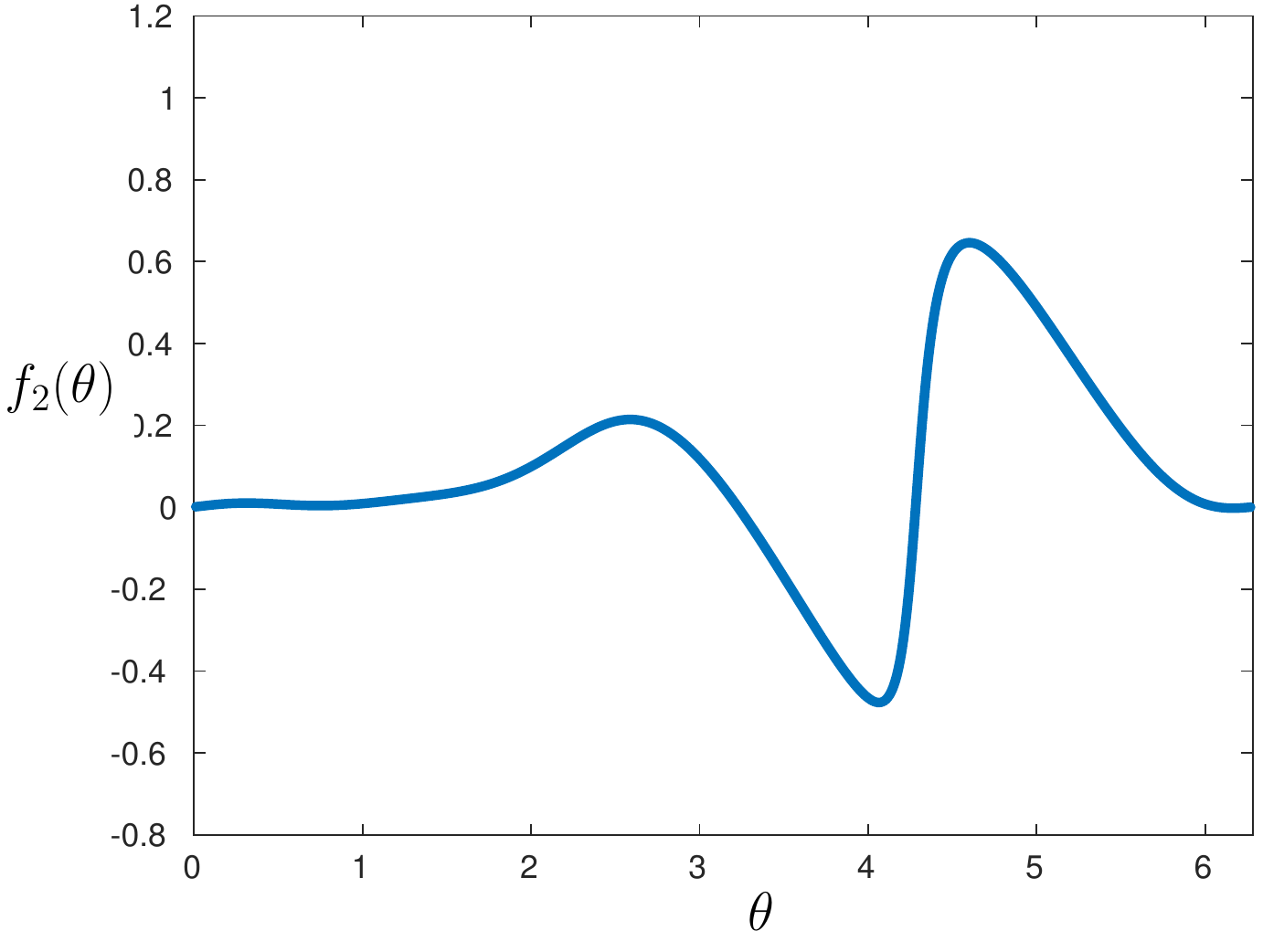}
\end{center}
\caption{Response function $f_2(\theta)$ which characterizes the phase response of a Hodgkin-Huxley neuron to a pulse with $u_{2max}$ corresponding to a current density of $10 \mu A/cm^2$, $u_{2min} = -u_{2max}/3$, and $p = 0.5$ ms.
\label{f_type2_10}}
\end{figure}

It will again be useful to consider the map which takes the phase of a neuron to its phase at a time $\tau$ later.  To formulate this map, we need the response curves for each type of pulse: the response curve $f(\theta)$ for the pulse with $u_{max}$ corresponding to $20 \mu A/cm^2$ was already shown in Figure~\ref{f_type2}; the response curve $f_2(\theta)$ for the pulse with $u_{max}$ corresponding to $10 \mu A/cm^2$ is shown in Figure~\ref{f_type2_10}.  To find this map, suppose that we start with $\theta(0^+) = 0$, immediately after the start of a pulse, where we assume that we have already accounted for the effect of the pulse according to the function $f(\theta)$. The next pulse, of different type, comes at time $\tau_2$.  Up until time $\tau_2$, the phase evolves according to $\dot{\theta} = \omega$; therefore,
\begin{equation}
\theta(\tau_2^-) = \theta_0 + \omega \tau_2.
\end{equation}
Treating the change in phase due to the next pulse as occurring instantaneously, we have
\begin{equation}
\theta(\tau_2^+) = \theta_0 + \omega \tau_2 + f_2(\theta_0 + \omega \tau_2).
\end{equation}
The system then evolves for a time $\tau - \tau_2$ without stimulus, giving
\begin{eqnarray*}
\theta(\tau^-) &=& \theta_0 + \omega \tau_2 + \omega (\tau - \tau_2) + f_2 (\theta_0 + \omega \tau_2) \\
&=& \theta_0 + \omega \tau + f_2 (\theta_0 + \omega \tau_2).
\end{eqnarray*}
At time $\tau$, we have another pulse of the type that started at $t=0$, so
\[
\theta(\tau^+) = \theta_0 + \omega \tau + f_2(\theta_0 + \omega \tau_2) + f(\theta_0 + \omega \tau + f_2(\theta_0 + \omega \tau_2)).
\]
Continuing in this fashion, we obtain
\begin{eqnarray*}
\theta(\tau + \tau_2^-) = \theta_0 &+& \omega (\tau + \tau_2) + f_2(\theta_0 + \omega \tau_2) \\
&+& f(\theta_0 + \omega \tau + f_2(\theta_0 + \omega \tau_2)),
\end{eqnarray*}
\begin{eqnarray*}
\theta(\tau + \tau_2^+) &=& \theta_0 + \omega (\tau + \tau_2) + f_2(\theta_0 + \omega \tau_2) \\
&& + f(\theta_0 + \omega \tau + f_2(\theta_0 + \omega \tau_2)) \\
&& + f_2(\theta_0 + \omega (\tau + \tau_2) + f_2(\theta_0 + \omega \tau_2) \\
&& \;\;\;\;\;\;\;\;\;\;\; + f(\theta_0 + \omega \tau + f_2(\theta_0 + \omega \tau_2)))
\end{eqnarray*}
\begin{eqnarray*}
\theta(2 \tau^-) &=& \theta_0 + 2 \omega \tau + f_2(\theta_0 + \omega \tau_2) \\
&& + f(\theta_0 + \omega \tau + f_2(\theta_0 + \omega \tau_2)) \\
&& + f_2(\theta_0 + \omega (\tau + \tau_2) + f_2(\theta_0 + \omega \tau_2) \\
&& \;\;\;\;\;\;\;\;\;\;\; + f(\theta_0 + \omega \tau + f_2(\theta_0 + \omega \tau_2)))
\end{eqnarray*}
\[
\theta(2 \tau^+) = \theta(2 \tau^-) + f(\theta(2 \tau^-)).
\]

A useful formulation is to let
\begin{equation}
G(s) = s + \omega \tau + f_2(s + \omega \tau_2) + f(s + \omega \tau + f_2(s + \omega \tau_2)),
\end{equation}
which gives
\begin{equation}
\theta(n \tau^+) = G^{(n)} (\theta_0).
\end{equation}
Alternatively, we can view this as a composition of two maps:
\[
\theta(0^+) = \theta_0,
\]
\[
\theta(\tau_2^+) = \theta_0 + \omega \tau_2 + f_2 (\theta_0 + \omega \tau_2) \equiv h_2(\theta_0),
\]
\begin{eqnarray*}
\theta(\tau^+) = \theta(\tau_2^+) + \omega (\tau-\tau_2) \!\! &+& \!\! f(\theta(\tau_2^+) + \omega (\tau- \tau_2)) \\
\equiv h_1(\theta(\tau_2^+)) &=& h_1(h_2(\theta_0)) = G(\theta_0).
\end{eqnarray*}
Note that we have written $G$, which is a map over the time interval $\tau$, as the composition of two maps $h_1$ and $h_2$, that is, 
\[
G = h_1 \circ h_2.
\]
Similar to before, we will look for fixed points of $G^{(n)}$, that is, solutions to $\theta^* = G^{(n)}(\theta^*)$.  If 
\begin{equation}
\left| \left. \frac{d}{d \theta} \right|_{\theta = \theta^*} (G^{(n)} (\theta)) \right|<1,
\end{equation}
then the fixed point of $G^{(n)}$ is stable.  Note that the relationship between fixed points of $G^{(n)}$ and clusters is more subtle for pulses with alternating properties than the relationship between fixed points of $g^{(n)}$ and $n$-clusters for identical pulses, because each $\tau$-interval for the alternating case contains two pulses.  This will be illustrated in the following examples.

\begin{figure}[tb]
\begin{center}
\leavevmode
\epsfxsize=2.2in
\epsfbox{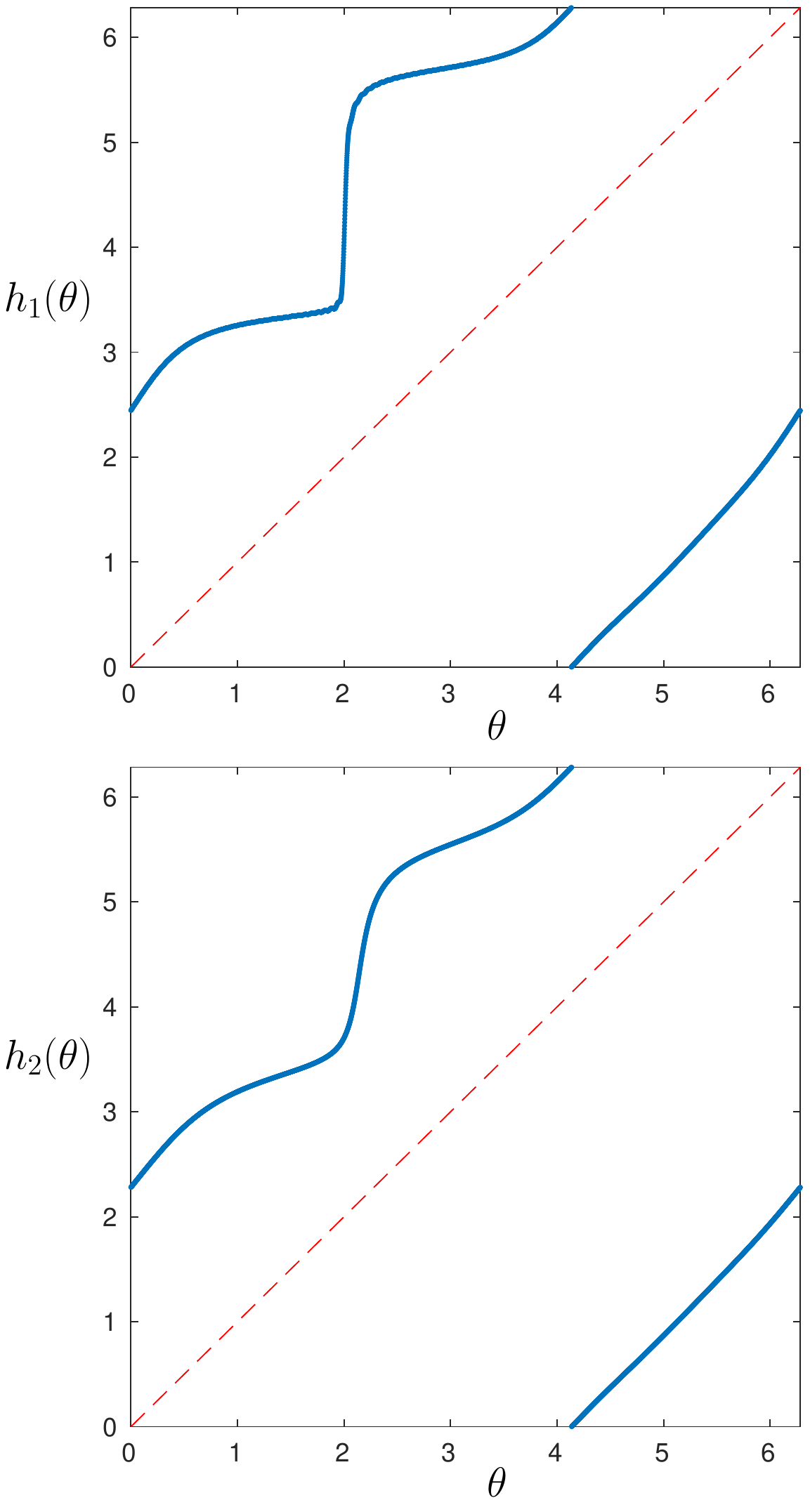}
\end{center}
\caption{Functions $h_1(\theta)$ for $u_{max}$ corresponding to a current density of $20 \mu A/cm^2$ and $h_2(\theta)$ for $u_{2max}$ corresponding to a current density of $10 \mu A/cm^2$, for $\tau = 10$ ms and $\tau_2 = \tau/2$.
\label{h1_h2_f100}}
\end{figure}
\begin{figure}[tb!]
\begin{center}
\leavevmode
\epsfxsize=2.2in
\epsfbox{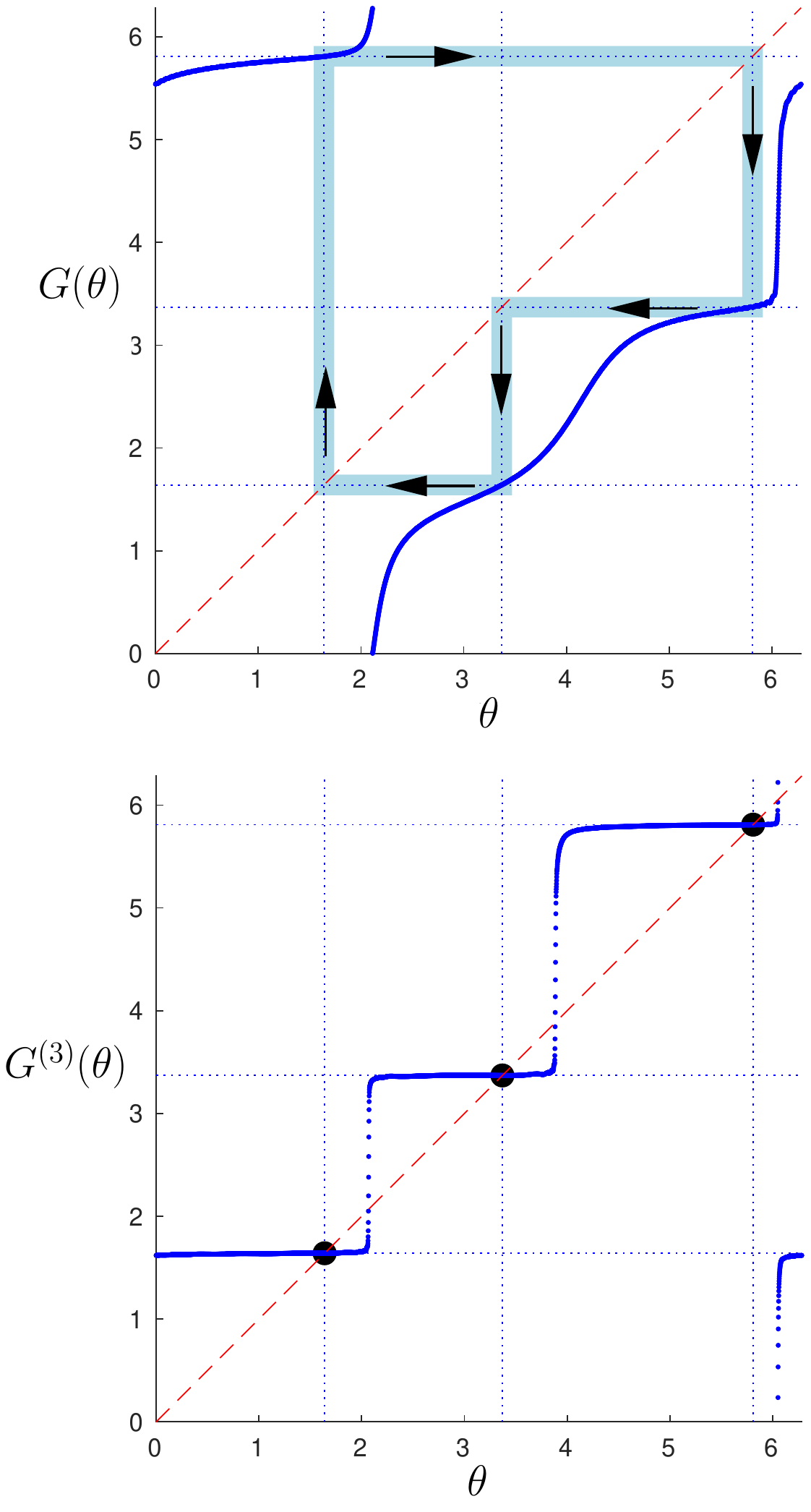}
\end{center}
\caption{Functions $G(\theta)$ and $G^{(3)}(\theta)$ for pulses with alternating properties with $u_{max}$ corresponding to a current density of $20 \mu A/cm^2$ and $u_{2max}$ corresponding to a current density of $10 \mu A / cm^2$, and $\tau = 10$ ms, $\tau_2 = \tau/2$.
\label{alternating_20_10_f100}}
\end{figure}

Figure~\ref{h1_h2_f100} shows $h_1(\theta)$ for $u_{max}$ corresponding to a current density of $20 \mu A/cm^2$ and $h_2(\theta)$ for $u_{2max}$ corresponding to a current density of  $10 \mu A/cm^2$, for $\tau = 10$ ms and $\tau_2 = \tau/2$.  We notice that these functions are quite similar to each other.  
Next, we show $G(\theta) = h_1(h_2(\theta))$ and $G^{(3)}(\theta)$ in Figure~\ref{alternating_20_10_f100}.  We see that there is a stable period-3 orbit for $G$, corresponding to three stable fixed points for $G^{(3)}$.  This corresponds to a 3-cluster state, as expected from Figure~\ref{type2_phase_vs_freq_two_kicks} evaluated at 100 Hz.  Here, the stable fixed points of $G^{(3)}$ correspond to a stable period-3 orbit of $G$, which in turn corresponds to a 3-cluster state.  We note that Figure~\ref{alternating_20_10_f100} looks very similar to Figure~\ref{f100} for identical pulses with frequency 100 Hz; however, the sequence of pulses is different.  For Figure~\ref{alternating_20_10_f100}, there is a ``large'' pulse at $t=0$, a ``small'' pulse at $t = 5$ ms, another large pulse at $t = 10$ ms, another small pulse at $t=15$ ms, another large pulse at $t = 20$ ms, etc.  For Figure~\ref{f100}, there is a large pulse at $t = 0$, another large pulse at $t = 10$ ms, another large pulse at $t = 20$ ms, etc, with no small pulses.

\begin{figure}[tb]
\begin{center}
\leavevmode
\epsfxsize=2.2in
\epsfbox{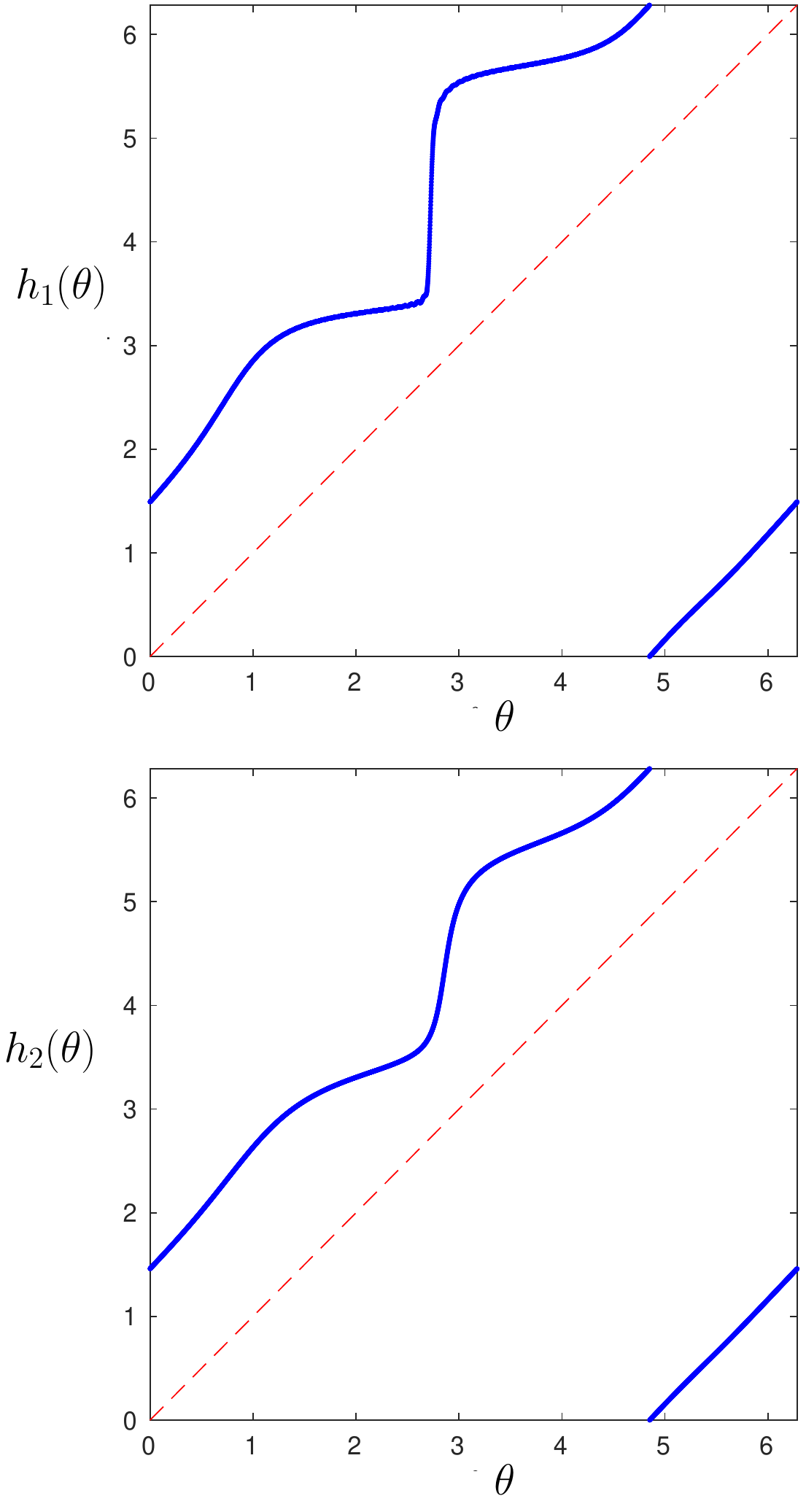}
\end{center}
\caption{Functions $h_1(\theta)$ for $u_{max}$ corresponding to a current density of $20 \mu A/cm^2$ and $h_2(\theta)$ for $u_{2max}$ corresponding to a current density of $10 \mu A/cm^2$, for $\tau = 6.67$ ms and $\tau_2 = \tau/2$.
\label{h1_h2_f150}}
\end{figure}

\begin{figure}[tb!]
\begin{center}
\leavevmode
\epsfxsize=2.2in
\epsfbox{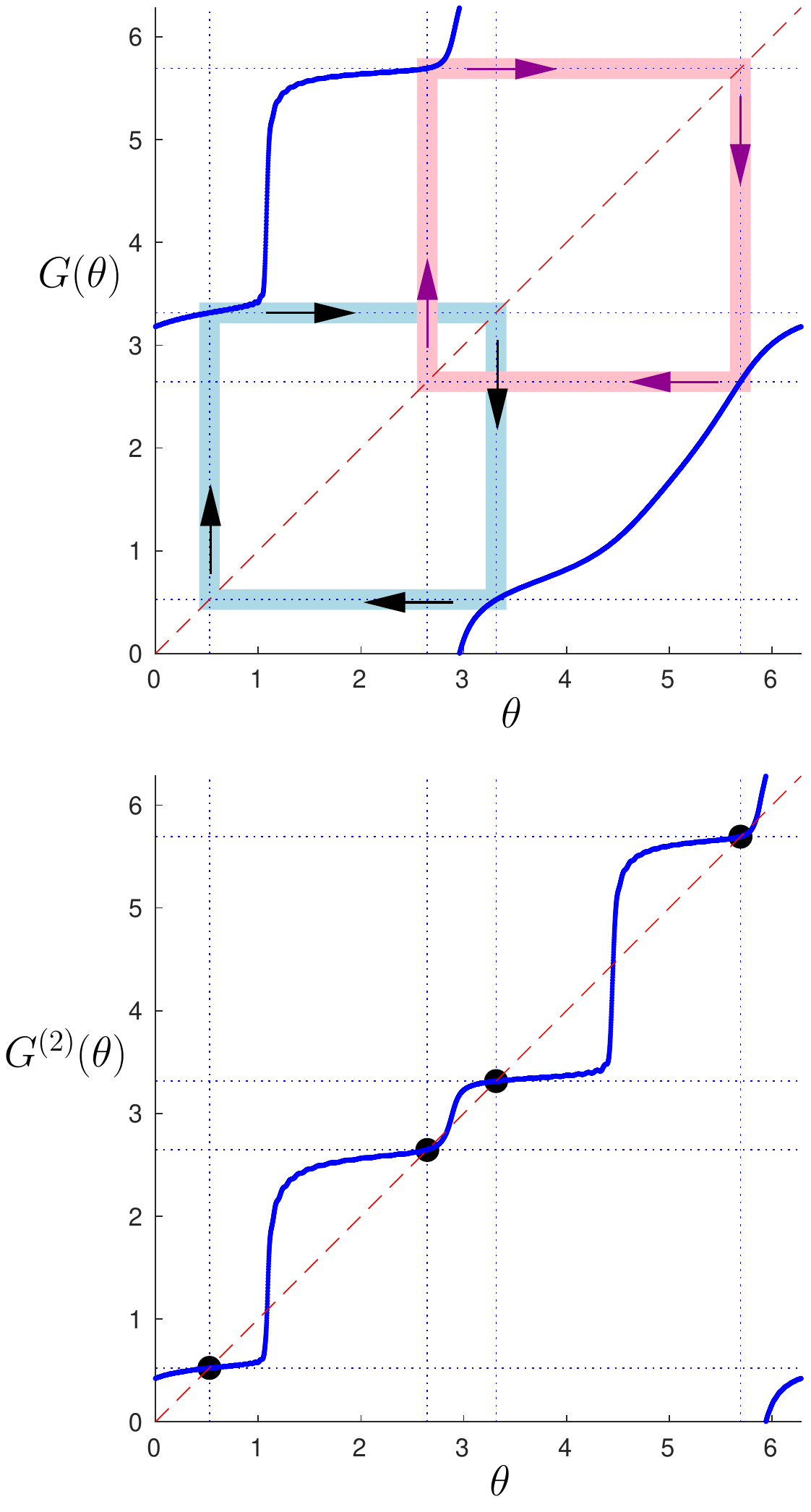}
\end{center}
\caption{Functions $G(\theta)$ and $G^{(2)}(\theta)$ for pulses with alternating properties with $u_{max}$ corresponding to a current density of $20 \mu A/cm^2$ and $u_{2max}$ corresponding to a current density of $10 \mu A / cm^2$, and $\tau = 6.67$ ms, $\tau_2 = \tau/2$.
\label{alternating_20_10}}
\end{figure}

\begin{figure}[h!]
\begin{center}
\leavevmode
\epsfxsize=3.5in
\epsfbox{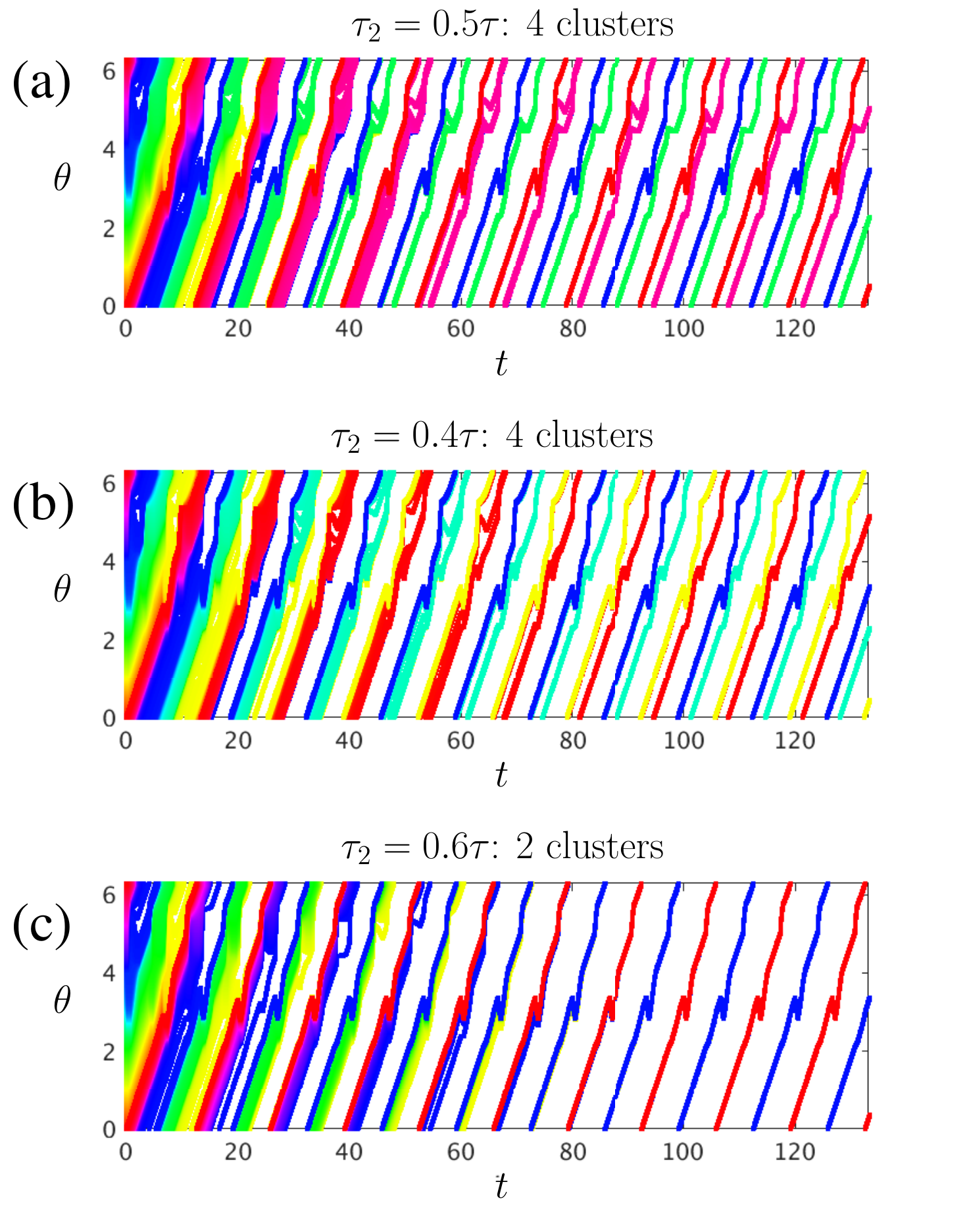}
\end{center}
\caption{Time series showing the phases of Hodgkin-Huxley neurons drawn from an initial uniform distribution with alternating pulses with $u_{max}$ corresponding to $20 \mu A/cm^2$ and $u_{2max}$ corresponding to $10 \mu A/cm^2$, for $\tau = 6.67$ ms and (a) $\tau_2 = 0.5 \tau$, (b) $\tau_2 = 0.4 \tau$, (c) $\tau_2 = 0.6 \tau$.  Four clusters form for (a) and (b), while only two clusters form for (c).
\label{phase_versus_time_two}}
\end{figure}

Figure~\ref{h1_h2_f150} shows $h_1(\theta)$ for $u_{max}$ corresponding to a current density of $20 \mu A/cm^2$ and $h_2(\theta)$ for $u_{2max}$ corresponding to a current density of  $10 \mu A/cm^2$, for $\tau = 6.67$ ms and $\tau_2 = \tau/2$.  As before, $h_1$ and $h_2$ are quite similar to each other.
Figure~\ref{alternating_20_10} shows $G(\theta) = h_1(h_2(\theta))$ and $G^{(2)}(\theta)$.  We see that there are four stable fixed points for $G^{(2)}$; these actually correspond to a 4-cluster state, as shown in Figure~\ref{type2_phase_vs_freq_two_kicks} evaluated at 150 Hz.  While at first it might seem surprising that stable fixed points for $G^{(2)}$ correspond to a 4-cluster state, we note that these results are similar to what we found for identical stimuli for a $300$ Hz stimulus (or, equivalently, for alternating pulses with $\tau = 6.67$ ms and  $u_{max} = u_{2max}$ corresponding to a current density of $20 \mu A/cm^2$, $\tau_2 = \tau/2$).  The proper comparison is that $G$ for alternating pulses with a stimulation frequency of $150$ Hz is similar to $g^{(2)}$ for identical pulses with a stimulation frquency of $300$ Hz, and $G^{(2)}$ for alternating pulses for a stimulation frequncy of $150$ Hz is similar to $g^{(4)}$ for identical pulses with a stimulation frequency of $300$ Hz.  These results show that we can obtain 4-cluster solutions for a population of oscillators with these alternating pulses; see Figure~\ref{phase_versus_time_two}(a).  


\begin{figure}[h!]
\begin{center}
\leavevmode
\epsfxsize=2.2in
\epsfbox{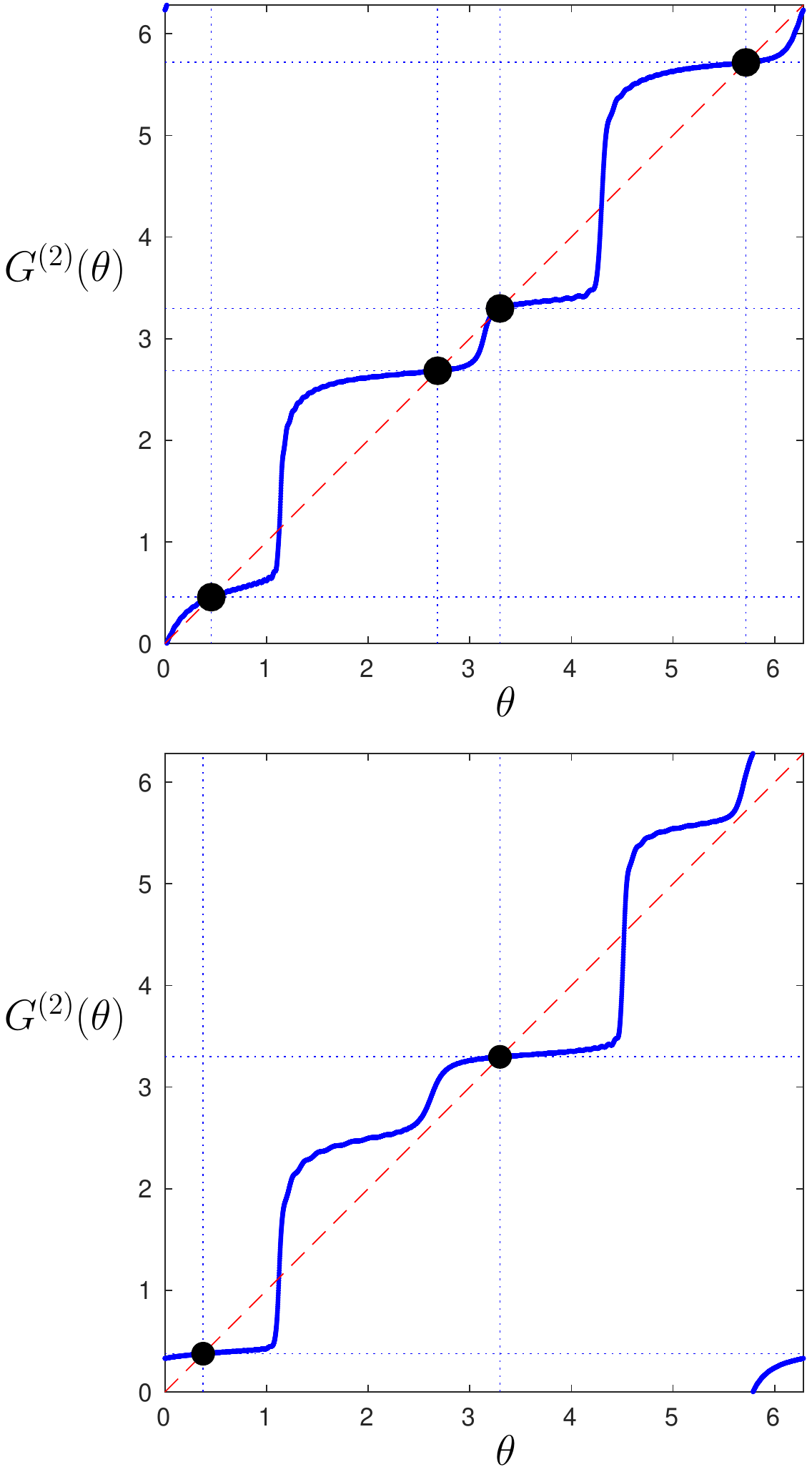}
\end{center}
\caption{Function $G^{(2)}(\theta)$ for alternating pulses with $u_{max}$ corresponding to $20 \mu A/cm^2$ and $u_{2max}$ corresponding to $10 \mu A/cm^2$, with $\tau = 6.67$ ms, and (left) $\tau_2 = 0.4 \tau$ and (right) $\tau_2 = 0.6 \tau$.
\label{alternating_20_10_tau0p4_tau0p6}}
\end{figure}

Our analytical formalism also allows one to consider alternating pulses for which $\tau_2 \neq \tau/2$.  For example, Figure~\ref{alternating_20_10_tau0p4_tau0p6} shows results for $u_{max}$ corresponding to $20 \mu A/cm^2$, $u_{2max}$ corresponding to $10 \mu A/cm^2$, and $\tau_2 = 0.4 \tau$ and $\tau_2 = 0.6 \tau$.  Interestingly, for $\tau_2 = 0.4 \tau$ there are four fixed points of the $G^{(2)}$ map, corresponding to a 4-cluster solution, but for $\tau_2 = 0.6 \tau$ there are only two fixed points of the $G^{(2)}$ map, corresponding to a 2-cluster solution.  Figures~\ref{phase_versus_time_two}(b) and (c) show the corresponding time series for these cases.   Comparing Figure~\ref{alternating_20_10_tau0p4_tau0p6} with the right panel of Figure~\ref{alternating_20_10}, we deduce that if $\tau_2$ is treated as a bifurcation parameter, there is a saddlenode bifurcation (this could also be called a tangent bifurcation for the $G^{(2)}$ map) for $\tau_2$ slightly larger than $0.5$.

\section{Conclusion} \label{section:conclusion}

Populations of neural oscillators subjected to periodic pulsatile stimuli 
can display interesting clustering behavior, in which subpopulations of the 
neurons are 
synchronized but the subpopulations are desynchronized with respect to each 
other.  The details of the clustering behavior depend on the frequency and 
amplitude of the stimuli in a complicated way.  Such clustering may be 
an important mechanism by which deep brain stimulation can lead to the 
alleviation of symptoms of Parkinson's disease and other disorders.

In this paper, we illustrated how the details of clustering for phase models 
of neurons subjected to periodic pulsatile inputs can be understood in terms 
of one-dimensional maps defined on the circle.  In particular, the analysis
allows one to predict the number of clusters, their stability properties, 
and their basins of attraction.  Moreover, we generalized our analysis to 
consider stimuli with alternating properties, which provide additional 
degrees of freedom in the design of DBS stimuli.

As part of our study, we found multiple ways to get the same type 
of clustering behavior, for example by using identical pulses or pulses 
with alternating properties, or from stimuli with different parameters such 
as stimulation frequency or the time spacing between pulses with alternating 
properties.  Such clustering occurs through the use of a single stimulation
electrode, unlike coordinated reset which requires multiple electrodes.  
We expect that the same clustering behavior can also be obtained for
different amplitudes of the pulses, 
cf.~\cite{wils15cluster}.  We believe that the analysis techniques used
in this paper can be useful for identifying a collection of stimuli
which give the same desireable clustering dynamics for a population of 
neurons, which will make it easier to find stimuli which are effective 
while minimizing the severity of side effects for DBS treatments.

Our analysis assumed certain properties of a neural population: all neurons 
are identical, they all receive the same input, they are uncoupled, and 
there is no noise.  For real neural populations, none of these assumptions
would be valid.  We also assumed that the phase models accurately capture
the dynamics of the neurons, which is only true for sufficiently small 
inputs; see, for example,~\cite{wils18}.  However, we believe that the 
results presented here form
an important baseline for the analysis of more realistic neural populations
stimulated by periodic pulses.  We note that the effect of noise on 
periodically forced neural populations has been considered in~\cite{wils15cluster}, 
which shows that for weak noise and long times, the number of neurons in each 
cluster is roughly the same.  We expect that similar results will hold
for neurons in the presence of weak noise subjected to alternating stimuli.

Our hope is that the techniques in this paper will help to guide the design of
stimuli for the treatment of Parkinson's disease and other disorders.  We 
believe that the use of pulses with alternating properties is particularly 
worthy of further investigation, since it represents a larger class of 
stimuli than has been considered in previous studies.

\section{Acknowledgements}
This research grew out of the Research Mentorship Program at the University of California, Santa Barbara during summer 2018.  We thank Dr.~Lina Kim for providing the opportunity for Daniel and Jacob to conduct this research as high school students, and for Tim Matchen for guidance on the project. 

\section*{Appendix: Neuron Models}

In this Appendix, we give details of the neural models used in the main text.

\vskip 0.1in
\noindent
{\it Thalamic neuron model}
\vskip 0.1in

The full thalamic neuron model is given by:

\begin{eqnarray*}
\dot V&=&\frac{-I_L-I_{Na}-I_K-I_T+I_b}{C_m}+u(t),\\
\dot h&=&\frac{h_{\infty}-h}{\tau_h},\\
\dot r&=&\frac{r_{\infty}-r}{\tau_r},
\end{eqnarray*}
where
\begin{eqnarray*}
    h_\infty &=& 1/(1+\exp((V+41)/4)),\\
    r_\infty &=& 1/(1+\exp((V+84)/4)),\\
    \alpha_h &=& 0.128\exp(-(V+46)/18),\\
    \beta_h &=& 4/(1+\exp(-(V+23)/5)),
    \end{eqnarray*}
    \begin{eqnarray*}
    \tau_h &=& 1/(\alpha_h+\beta_h),\\
    \tau_r &=& (28+\exp(-(V+25)/10.5)),\\
\end{eqnarray*}
\begin{eqnarray*}
    m_\infty &=& 1/(1+\exp(-(V+37)/7)),\\
    p_\infty &=& 1/(1+\exp(-(V+60)/6.2)),\\
\end{eqnarray*}
\begin{eqnarray*}
    I_L&=&g_L(v-e_L),\\
     I_{Na}&=&g_{Na}({m_\infty}^3)h(v-e_{Na}),\\
    I_K&=&g_K((0.75(1-h))^4)(v-e_K),\\
    I_T&=&g_T(p_\infty^2)r(v-e_T).
\end{eqnarray*}
The parameters for this model are
\[
C_m = 1 \; \mu F/cm^2 \;,\;   g_L = 0.05 \; mS/cm^2 \;,\;  e_L = -70 \; mV \;,
\]
\[
g_{Na} = 3 \; mS/cm^2 \;,\;   e_{Na} = 50 \; mV, g_K = 5 \; mS/cm^2 \;,
\]
\[
e_K = -90 \; mV \;,\;   g_T = 5 \; mS/cm^2 \;,\;   e_T = 0 \; mV \;, 
\]
\[
I_b = 5 \; \mu A/cm^2.
\]

\vskip 0.1in
\noindent
{\it Hodgkin-Huxley neuron model}
\vskip 0.1in

The full Hodgkin-Huxley model is given by:

\begin{eqnarray*}
\dot{V}&=&(I_b -\bar{g}_{Na}h(V-V_{Na})m^3-\bar{g}_K(V-V_K)n^4\\
&& -\bar{g}_L(V-V_L))/c + u(t)  \; , \\
\dot{m}&=& a_m(V)(1-m)-b_m(V)m  \; , \\
\dot{h}&=&a_h(V)(1-h)-b_h(V)h  \; , \\
\dot{n}&=&a_n(V)(1-n)-b_n(V)n  \; , 
\end{eqnarray*}
where
\begin{eqnarray*}
a_m(V) &=&  0.1(V+40)/(1-\exp(-(V+40)/10))  \; , \\
b_m(V) &=&  4\exp(-(V+65)/18)  \; , \\
a_h(V) &=&  0.07\exp(-(V+65)/20)  \; , \\
b_h(V) &=&  1/(1+\exp(-(V+35)/10))  \; , \\
a_n(V) &=&  0.01(V+55)/(1-\exp(-(V+55)/10))  \; , \\
b_n(V) &=&  0.125\exp(-(V+65)/80) \; , 
\label{eq_HH}
\end{eqnarray*}
The parameters for this model are
\[
V_{Na}=50 \; mV \;,\; V_K=-77 \; mV \;,\; V_L=-54.4 \; mV \;, \\
\]
\[
\bar{g}_{Na}=120 \; mS/cm^2  \; , \; \bar{g}_K=36 \; mS/cm^2 \;, 
\]
\[
bar{g}_L=0.3 \; mS/cm^2 \;,\; I_b=10 \; \mu A/cm^2 \;,
\]
\[
c=1 \; \mu F/cm^2 .
\]

\vspace{-0.2in}


\bibliographystyle{spmpsci}      

\bibliography{ms}



\end{document}